\documentclass[12pt,oneside]{amsart}

\usepackage{amssymb,mathrsfs,mathcomp,amsfonts,longtable, hhline}
\usepackage{url} 
\usepackage[matrix,arrow,curve]{xy}

\makeatletter \@addtoreset{equation}{section} \makeatother

\newcommand{\muu}{{\boldsymbol{\mu}}}

\newcommand{\ct}{\operatorname{ct}}

\newcommand{\Pic}{\operatorname{Pic}}
\newcommand{\DP}{\operatorname{DP}_{5}^{\operatorname{A}_4}}

\newcommand{\Sing}{\operatorname{Sing}}

\newcommand{\Bs}{\operatorname{Bs}}
\newcommand{\Cl}{\operatorname{Cl}}
\newcommand{\rk}{\operatorname{rk}}
\newcommand{\g}{\operatorname{g}}

\newcommand{\qW}{\operatorname{qW}}
\newcommand{\qQ}{\operatorname{q\QQ}}

\newcommand{\B}{{\mathbf B}}
\newcommand{\CC}{\mathbb{C}}
\newcommand{\FF}{\mathbb{F}}
\newcommand{\QQ}{\mathbb{Q}}

\newcommand{\ZZ}{\mathbb{Z}}
\newcommand{\PP}{\mathbb{P}}
\newcommand{\OOO}{{\mathscr{O}}} 
\newcommand{\HHH}{{\mathscr{H}}} 
\newcommand{\MMM}{{\mathscr{M}}} 
 
\newcommand{\SSS}{{\mathscr{S}}}
\renewcommand{\emptyset}{\varnothing}

\newcommand{\qq}{\mathbin{\sim_{\scriptscriptstyle{\QQ}}}}
\newcommand{\notqq}{\mathbin{{\not\sim}_{\scriptscriptstyle{\QQ}}}}

\newcounter{NO}

\newcommand{\comment}[1]{}

\newcommand{\xref}[1]{{\rm \ref{#1}}}

\newtheorem{theorem}{Theorem}

\numberwithin{theorem}{section}
\numberwithin{equation}{theorem}

\newtheorem{mtheorem}[theorem]{}
\newtheorem{stheorem}[equation]{}

\theoremstyle{definition}
\newtheorem{case}[theorem]{}
\newtheorem{scase}[equation]{}

\newtheorem*{remark*}{Remark}

\newcounter{NN}\numberwithin{NN}{section}
\renewcommand{\theNN}{\arabic{NN}${}^o$}

\def\nr{\refstepcounter{NN}{\theNN}}%

\title{$\QQ$-Fano threefolds of large Fano index, II}

\author{Yuri Prokhorov}
\thanks{
The author's work partially supported by 
RFBR grants No. 11-01-00336-a, 11-01-92613-KO\_a, the grant of
Leading Scientific Schools No. 4713.2010.1 and 
AG Laboratory SU-HSE, RF government 
grant ag. 11.G34.31.0023.
}

\address{Department 
of Algebra, Faculty of Mathematics, Moscow State
University, Moscow 117234, Russia
\newline\indent
Laboratory of Algebraic Geometry, SU-HSE, 
7 Vavilova Str., Moscow 117312, Russia
}
\email{prokhoro@gmail.com}
\begin{document}
\begin{abstract}
We classify Q-Fano threefolds of Fano index $> 2$ and 
big degree.
\end{abstract}

\maketitle

\section{Introduction}\label{sect-intr}
This work is a sequel to our previous papers \cite{Prokhorov-2007-Qe},
\cite{Prokhorov2008a}. 
Recall that a three-dimensional projective variety $X$ 
is called \textit{$\QQ$-Fano threefold} if it has only terminal 
$\QQ$-factorial singularities, $\Pic(X)\simeq \ZZ$,
and its anticanonical divisor $-K_X$ is ample.
In this situation we define \textit{Fano-Weil} and \textit{$\QQ$-Fano indices}
as follows:
\[
\begin{array}{lll}
\qW(X)&:=&\max \{ q\in\ZZ \mid -K_X\sim q A,\ \text{$A$ is a Weil divisor}\},
\\[5pt]
\qQ(X)&:=&\max \{ q\in \ZZ \mid -K_X\qq q A,\ \text{$A$ is a Weil divisor}\},
\end{array}
\]
where $\sim$ (resp. $\qq$) denotes 
linear (resp. numerical) equivalence.
Clearly, $\qW(X)$ divides $\qQ(X)$.
Another important invariant of a Fano variety $X$ is 
the \textit{genus} $\g(X):=\dim |-K_X|-1$.
It is known that \refstepcounter{theorem}
\begin{equation}
\label{equation-index}
\qW(X),\ \qQ(X)\in \{1,\dots, 11, 13,17,19\}
\end{equation}
(see \cite{Suzuki-2004}, \cite[Lemma 3.3]{Prokhorov2008a}).
Moreover, in \cite{Prokhorov2008a} we proved the following result.

\begin{mtheorem}{\bf Theorem.}
\label{theorem-main-10-19}
Let $X$ be a $\QQ$-Fano threefold with $q:=\qQ(X)\ge 9$. Then $q=\qW(X)$ and $\Cl(X)\simeq \ZZ$.
\begin{enumerate}
\item\label{theorem-main-q=19}
If $q=19$, then $X\simeq \PP(3,4,5,7)$.
\item\label{theorem-main-q=17}
If $q=17$, then $X\simeq \PP(2,3,5,7)$.

\item\label{theorem-main-q=13}
If $q=13$ and $\g(X) >4$, then $X\simeq \PP(1,3,4,5)$.

\item\label{theorem-main-q=11}
If $q=11$ and $\g(X) >10$, then $X\simeq \PP(1,2,3,5)$.

\item\label{theorem-main-q=10}
$q\neq 10$.
\end{enumerate}
\end{mtheorem}
In this paper we study $\QQ$-Fano threefolds with
$3\le \qQ(X)\le 9$ and ``sufficiently big'' genus.
Our main result is the following.

\begin{mtheorem}{\bf Theorem.}
\label{theorem-main}
Let $X$ be a $\QQ$-Fano threefold and let $q:=\qQ(X)$.

\begin{enumerate}
\item\label{theorem-main-q=9}
If $q=9$ and $\g(X)>4$, then $\g(X)=18$. In this case,
$X\simeq X_6 \subset \PP(1,2,3,4,5)$ and 
the equation of $X_6$, in suitable coordinates, can be written as
$x_1x_5+x_2x_4+x_3^2=0$.

\item\label{theorem-main-q=8}
If $q=8$ and $\g(X)>10$, then either
\begin{enumerate}
\item
$X\simeq X_6\subset \PP(1,2,3^2,5)$ and the equation of
$X_6$, in suitable coordinates, can be written as
$x_1x_5+x_3x_3'+x_2^3=0$, or
\item
$X\simeq X_{10}\subset\PP(1,2,3,5,7)$.
\end{enumerate}

\item\label{theorem-main-q=7}
If $q=7$ and $\g(X)>17$, then 
$X\simeq \PP(1^2,2,3)$.

\item\label{theorem-main-q=6}
If $q=6$ and $\g(X)>15$, then
$X\simeq X_6\subset \PP(1^2,2,3,5)$ and the equation of
$X_6$, in suitable coordinates, can be written as
$x_1'x_5+x_3^2+x_2^3+\lambda x_1^4x_2+\mu x_1^6=0$.

\item\label{theorem-main-q=5}
If $q=5$ and $\g(X)>18$, then
$X\simeq \PP(1^3,2)$ or $X_4\subset \PP(1^2,2^2,3)$.

\item\label{theorem-main-q=4}
If $q=4$ and $\g(X)>21$, then
$X\simeq \PP^3$ or $X_4\subset \PP(1^3,2,3)$.

\item\label{theorem-main-q=3}
If $q=3$ and $\g(X)>20$, then
$X\simeq X_2\subset \PP^4$ or $X_3\subset \PP(1^4,2)$.
\end{enumerate}
\end{mtheorem}
The following consequence of our theorem was inspired by 
discussions with J\"urgen Hausen.

\begin{mtheorem} {\bf Corollary.}
\label{corollary-main}
Let $\mathcal X_q$ be the class of all $\QQ$-Fano threefolds $X$ with $\qQ(X)=q\ge 3$
and let $X_q\in \mathcal X_q$ be a $\QQ$-Fano threefold of maximal 
genus $\g(X)$.
Then $X_q$ admits an effective action of a torus of dimension $\ge 2$. 
\end{mtheorem}
In proofs we use a computer computation based on 
the orbifold Riemann-Roch theorem and 
Bogomolov-Miyaoka inequality. It gives us 
472 possible candidates of Hilbert series of $\QQ$-Fano threefolds with $\qQ(X)=\qW(X)\ge 3$.
The most essential part of the paper is to disprove some cases or
to show that some particular $\QQ$-Fano threefold is isomorphic to
a certain hypersurface in a weighted projective space. Here we use 
techniques which is a generalization 
of Fano-Iskovskikh ``double projection'' method. 
It was adopted to the singular case by Alexeev \cite{Alexeev-1994ge}
(see also \cite{Prokhorov-2007-Qe},
\cite{Prokhorov2008a}). 
The method can be applied in many other cases.
For example we prove the following.

\begin{mtheorem}{\bf Theorem.}
\label{main-hypersurfaces}
Let $X$ be a $\QQ$-Fano threefold.
\begin{enumerate}
 \item 
If $\qQ(X)=7$, $\g(X)=17$, and $-K_X^3=7^3/10$, then
$X\simeq X_6\subset \PP(1,2^2,3,5)$.
 \item 
 The case $\qQ(X)=8$, $\g(X)=10$, and $-K_X^3=8^3\cdot 2/45$ does not occur.
\end{enumerate}
\end{mtheorem}
Note however that typically our method works for Fano threefolds having 
sufficiently many anticanonical sections.
For another approaches 
we refer to \cite{Reid2003}, \cite{Brown2010}. 

The progress in the classification of $\QQ$-Fano threefolds of 
Fano index $\ge 3$ is summarized in the following table (for simplicity we assume 
here that the group $\Cl(X)$ is torsion free):
\begin{center}
\begin{tabular}{p{120pt}cccccccccccc}
$\qW(X)$ & 3 & 4 & 5 & 6 & 7 & 8 & 9 & 10&11 & 13 & 17 & 19
\tabularnewline[5pt]
\hline 
\tabularnewline[-5pt]
\# candidates & 231 & 124 & 63 & 11 & 23 & 10 & 2&1 & 3 & 2 & 1 & 1
\tabularnewline[5pt]
of which do not occur & 12 & 11 & 11 & 2 & 7 & 5 &0& 1 & 0 & 0 & 0 & 0
\tabularnewline[5pt]
of which exist & 13 & 10 & 9 & 1 & 5 & 3 & 2 & 0&3 &2 & 1 & 1
\tabularnewline[5pt]
of which exist and \newline completely described & 4 & 3 & 3 & 1 & 2 & 2&1 & 0 & 1 & 1 & 1 & 1
\end{tabular}
\end{center}
\par \medskip\noindent
In this table the second row shows the number of admissible baskets 
of corresponding Fano index according to 
\cite{GRD}, \cite{Brown-Suzuki-2007j}, \cite{Prokhorov-2010-qfano}.
In the third row we indicate the number of candidates which 
do not occur as $\QQ$-Fano threefolds. This is according to this paper
and \cite{Prokhorov-2007-Qe}.
The fourth row shows the number baskets for which we know examples of 
$\QQ$-Fano threefolds. Here all known Fano threefolds are of codimension 
$\le 3$ in some weighted projective space \cite{GRD}, \cite{Brown-Suzuki-2007j}.
The last row shows the number baskets which are completely investigated:
we know the explicit description of corresponding $\QQ$-Fano threefolds
according to \cite{Sano-1996}, this paper, and \cite{Prokhorov2008a}.

Codimension $4$ \ $\QQ$-Fano threefolds of 
were investigated in \cite{Brown2010}
under an assumption that the Fano index equals to 1.
The authors proved that in many cases a codimension $4$ basket 
can be realized by at least two different deformation families of 
$\QQ$-Fano threefolds.
There is a hope to adopt this method to 
$\QQ$-Fanos of higher Fano index. 
Note however that it follows from our main theorem that the following 
codimension $4$ candidates do not exist:
\begin{itemize}
 \item 
$\qW(X)=7$, $\B(X)=(2,2,8)$, $A^3=1/8$, 
\cite[no. 41482]{GRD},
 \item 
$\qW(X)=6$, $\B(X)=(7)$, $A^3=2/7$, 
\cite[no. 41468]{GRD},
 \item 
$\qW(X)=4$, $\B(X)=(5)$, $A^3=4/5$, 
\cite[no. 41356]{GRD},
 \item 
$\qW(X)=3$, $\B(X)=(5)$, $A^3=8/5$, 
\cite[no. 41225]{GRD}.
\end{itemize}

\textbf{Acknowledgments.}
This paper grew out of a talk given at the
``Higher Dimensional Birational Geometry'' workshop
at the University of Warwick
during my visit in June, 2010.
I would like to thank 
the Warwick Institute of Advanced Study for the support 
and Professor Miles Reid for the invitation, hospitality and his interest in current 
computation.

\section{Preliminaries}
\label{sect-2}
\begin{case} {\bf Notation.}
We work over the complex number field $\CC$ throughout.

\begin{itemize}
\item[] $\B(X)$
denotes the basket of a terminal threefold $X$,
\item[] $\Cl(X)$ denotes the Weil divisor class group,
\item[] $\rho(X):=\rk \Pic(X)$,
\item[] $\g(X):=\dim |-K_X|-1$, the \textit{genus} of a Fano threefold $X$,
\item[] $\PP(a_1,\dots, a_n)$ is the weighted projective space,
\item[] $X_d\subset \PP(a_1,\dots, a_n)$ is a hypersurface of weighted degree $d$.
\end{itemize} 
\end{case} 
 
\begin{case}
\label{computer} 
We will use systematically a computer computation based on
the orbifold Riemann-Rich theorem \cite{Reid-YPG1987}
and Bogomolov-Miyaoka inequality (see \cite{Kawamata-1992bF}).
The algorithm is explained in \cite{Brown-Suzuki-2007j} and \cite[Proof of Lemma 3.5]{Prokhorov2008a}.
It allows us to find all possible baskets and Hilbert series of $\QQ$-Fano threefolds
with $\qQ(X)=\qW(X)\ge 3$. The corresponding computer programs can be found 
in the Graded Ring Database \cite{GRD} (MAGMA code) or on author's homepage \cite{Prokhorov-2010-qfano}
(PARI code). Note that we may assume that $-K_X^3\le 125/2$ and the equality holds only for 
$X\simeq \PP(1^3,2)$ \cite{Prokhorov-2007-Qe}.
\end{case}

To apply \ref{computer} we typically have to show that 
the group $\Cl(X)$ is torsion free.
This can be done by using the following.

\begin{mtheorem} {\bf Lemma (cf. \cite[Prop. 5.3]{Prokhorov-2007-Qe}).}
\label{lemma-torsion}
Let $X$ be a $\QQ$-Fano threefold with $\qQ(X)>3$.
Assume that $\Cl(X)$ contains an element $\xi$ of order $n>1$.
Then there exists a $\QQ$-Fano threefold $Y$ such that $\qQ(X)$ divides $\qQ(Y)$
and $-K_Y^3\ge -nK_X^3$. 
Moreover,
\[
\g(Y)\ge n
(\g(X)-1)-3.
\]
\end{mtheorem}
\begin{proof}
Similar to the proof of \cite[Prop. 5.3]{Prokhorov-2007-Qe}.
\end{proof}

\begin{case}
\label{construction-notation}
\textbf{Construction.}
We use notation of \cite{Prokhorov2008a}.
Let $X$ be a $\QQ$-Fano threefold and let $A$ be an
ample Weil divisor whose class generates the group $\Cl(X)/{\qq}$. 
Thus we have $-K_X \qq qA$.

In the construction below we follow \cite{Alexeev-1994ge}.
All the facts are explained in \cite{Prokhorov2008a}, so basically we omit proofs.
Let $\MMM$ be a mobile linear system 
without fixed components and let $c:=\ct(X,\MMM)$ be the canonical threshold of $(X,\MMM)$.
So the pair $(X,c\MMM)$ is canonical but not terminal.
Assume that $-(K_X+c\MMM)$ is ample.
Let $f\colon \tilde X\to X$ be a $K+c\MMM$-crepant blowup in the Mori category so that 
$\tilde X$ has only terminal $\QQ$-factorial singularities, $\rho(\tilde X/X)=1$, and 
\begin{equation}
\label{equation-Sarkisov-discr}
K_{\tilde X}+c\tilde \MMM=f^*(K_X+c\MMM).
\end{equation} 
Here the exceptional locus $E\subset \tilde X$ is an irreducible divisor.

As in \cite{Alexeev-1994ge}, run $K+c\MMM$-MMP on $\tilde X$. 
We get the following diagram \textup(Sarkisov link of type I or II\textup)
\begin{equation}
\label{eq3.1}
\begin{gathered}
\xymatrix{
&\tilde X\ar@{-->}[r]\ar[dl]_{f}&\bar X\ar[dr]^{\bar f}&
\\
X&&&\hat X
} 
\end{gathered}
\end{equation}
where the varieties~$\tilde X$ and~$\bar X$ have only $\QQ$-factorial terminal
singularities, $\rho(\tilde X)=\rho (\bar X)=2$, 
$f$ is a Mori extremal divisorial contraction,
$\tilde X \dashrightarrow \bar X$ is a sequence of 
log flips, and $\bar f$ is a (not necessarily birational) Mori extremal contraction.
In particular, $\rho(\hat X)=1$. Let $\Theta$ be a Weil divisor on $\hat X$ 
whose class generates $\Cl(\hat X)/\qq$. 
 
If the linear system $|kA|$ is not empty,
everywhere below we denote by $S_k\in |kA|$ a general member.
In all what follows, for a divisor (or a linear system) $D$ on $X$,
let $\tilde D$ and $\bar D$ denote 
strict birational transforms of $D$ on~$\tilde X$ and $\bar X$, respectively.

We can write
\begin{equation}
 \label{equation-beta}
f^*S_k=\tilde S_k+\beta_k E,
\end{equation}
\begin{equation}
\label{equation-discrepancies}
K_{\tilde X}=f^*K_X+\alpha E,\quad K_{\tilde X}+\tilde \MMM +a_0E=f^*(K_X+\MMM),
\end{equation}
where $\beta_k\ge 0$ and $\alpha, \, a_0>0$. 
If $\MMM=|kA|=|S_k|$, then $a_0=\beta_k-\alpha$.
Note that 
the divisor $-K_{\tilde X}\equiv c\tilde \MMM-f^*(K_X+c\MMM)$
is nef outside of a finite number of curves on $\tilde X$.
The same holds for $-K_{\bar X}$.

\begin{stheorem} {\bf Lemma ({\cite[4.2]{Prokhorov2008a}}).}
\label{lemma-c-42}
Let $P\in X$ be a point of index $r>1$. Assume that 
$\MMM\sim -mK_X$ near $P$, where $0<m<r$. Then $c\le 1/m$. 
\end{stheorem}
\end{case}
\begin{case}
\label{case-not-birational-del-Pezzo}
Assume that $\bar f$ is not birational.
Then $\hat X$ is either a smooth rational curve or 
a del Pezzo surface with at worst Du Val singularities and 
$\rho(\hat X)=1$ \cite{Mori-Prokhorov-2008}.
We also have $\bar f(\bar E)=\hat X$.
If moreover $\Cl(X)$ is torsion free, then $\hat X\simeq \PP^1$,
$\PP^2$, $\PP(1^2,2)$, $\PP(1,2,3)$, or $\DP$, where $\DP$ is a Du Val del Pezzo surface 
of degree $5$ with $\rho=1$ (see Proposition \ref{lemma-DuVal-5}).
The following obvious computation will be used frequently throughout the paper.
\[
\begin{array}{c|ccccc}
\hat X & \PP^1& \PP^2 & \PP(1^2,2) & \PP(1,2,3)& \DP
\\[3pt]
\hline
\\[-5pt]
\dim |\Theta|& 1&3&1&0&0
\end{array} 
\] 
\end{case}

\begin{case}
Assume that the contraction $\bar f$ is birational. Then $\hat X$ is a $\QQ$-Fano threefold.
In this case, denote by $\bar F$ the $\bar f$-exceptional divisor,
$\tilde F\subset \tilde X$ its proper transform, $F:=f(\tilde F)$, and 
$\hat q:=\qQ(\hat X)$. 
For a divisor $\bar D$ on $\bar X$, we
put $\hat D:=\bar f_*\bar D$.
Write $\hat S_k\qq s_k \Theta$ and $\hat E\qq e\Theta$, where
$s_k$ and $e$ are non-negative integers. 
One can see that $\bar E\neq \bar F$. Thus $e>0$.
Clearly, $s_k=0$ if and only if $\bar S_k= \bar F$.
In particular, $s_k>0$ if $S_k$ is moveable.
Similar to \eqref{equation-discrepancies} we write 
\begin{equation}
\label{equation-discrepancies-second}
K_{\bar X}=\bar f^*K_{\hat X}+b \bar F,\quad
\bar S_k=\bar f^*\hat S_k-\gamma_k \bar F,\quad
\bar E=\bar f^*\hat E-\delta \bar F,
\end{equation}
where $b$, $\delta$, $\gamma_k$ are non-negative rational numbers, and $b>0$.

\begin{stheorem}{\bf Lemma.}
\label{lemma-dk}
If $\bar f$ is birational and $\alpha<1$, then $\g(X)\le \g(\hat X)$. 
\end{stheorem}

\begin{proof}
Let $\HHH:=|-K_X|$. Similar to  \eqref{equation-discrepancies}
we can write
\[
K_{\tilde X}+\tilde \HHH +a_0'E=f^*(K_X+\HHH)\sim 0, 
\]
where $a_0'\ge -\alpha>-1$. Since $K_X+\HHH$ is Cartier,  we have $a_0'\in \ZZ$ and
$a_0'\ge 0$. Hence, $-K_{\hat X}\sim \hat \HHH +a_0'\hat E$
and so $\dim |-K_{\hat X}|\ge \dim \hat \HHH=\dim \HHH$.
\end{proof}

\begin{stheorem}{\bf Lemma ({\cite[4.12]{Prokhorov2008a}}).}
\label{lemma-torsion-ae}
Assume that $\Cl(X)$ is torsion free and the map $\bar f$ is birational.
Write $F\sim dA$.
Then $\Cl(\hat X)\simeq \ZZ\oplus \ZZ_{n}$, where $d=ne$. 
If additionally $s_1=0$ \textup(i.e., $\bar S_1$ is 
$\bar f$-exceptional\textup), then $\Cl(\hat X)$ is torsion free and $e=1$. 
\end{stheorem}
\end{case}

\begin{mtheorem}{\bf Lemma (\cite{Kawakita2005}).}
\label{lemma-discrepancies}
Let $X\ni P$ be a threefold terminal point of index $r>1$
with basket 
$\B(X,P)$ \textup(see \cite{Reid-YPG1987}\textup) and 
let $f: (\tilde X\supset E)\to (X\ni P)$ be a divisorial Mori extraction,
where $E$ is the exceptional divisor and $f(E)=P$. 
Write $K_{\tilde X}=f^*K_X+\alpha E$.
\begin{enumerate}
 \item 
If $X\ni P$ is a point of type other than $\mathrm{cA}/r$ and $r>2$, then 
$\alpha=1/r$.
 \item 
If $X\ni P$ is of type $\mathrm{cA}/r$
and $\B(X,P)$ consists of $n$ points of index $r$, then $\alpha=a/r$, where $n\equiv 0\mod a$.
\end{enumerate}
\end{mtheorem}
\begin{proof}
The assertion (i) immediately follows by \cite[Theorems 1.2 and 1.3]{Kawakita2005}.
Assume that $X\ni P$ is of type $\mathrm{cA}/r$ and $\alpha>1/r$. 
Again according to \cite[Theorems 1.2 and 1.3]{Kawakita2005}
the point $P\in X$ can be identified with
\[
\{x_1x_2+\phi(x_3^r,x_4)=0\}\subset \CC^4/\muu_r(1,-1,b,0), 
\]
so that $f$ is the weighted blowup with weights 
$(r_1/r, r_2/r, a/r,1)$, where
\begin{itemize}
\item
$a\equiv br_1\mod r$, and $r_1+r_2\equiv 0\mod ra$,
 \item
$(a- br_1)/r$ is coprime to $r_1$,
 \item
$\phi$ has weighted order $(r_1+r_2)/r$ with weights
$(a/r,1)$,
 \item
the monomial $x_3^{(r_1+r_2)/a}$ appears in $\phi$ with a nonzero coefficient.
\end{itemize}
Then $\alpha=a/r$.
In our case $x_4^n$ appears in $\phi$ with a nonzero coefficient.
Hence, $(r_1+r_2)/r= n$ and so $n\equiv 0\mod a$.
\end{proof}

\begin{mtheorem}{\bf Proposition.} 
\label{lemma-DuVal-5}
Let $Z$ be a del Pezzo surface of degree $5$ with 
only Du Val singularities.
The following are equivalent:
\begin{enumerate}
\item 
$\Pic(Z)\simeq \ZZ$,
\item
$\Cl(Z)\simeq \ZZ$,
\item 
$\Sing(Z)$ consists of one Du Val point of type $A_4$,
\item 
$-K_Z\qq 5 L$ for some Weil divisor $L$.
\end{enumerate}
Moreover, such $Z$ is unique up to isomorphism and isomorphic to 
a hypersurface of degree $6$ in $\PP(1,2,3,5)$.
\end{mtheorem}
\begin{proof}
Almost all assertions are well-known (see e.g. \cite[\S 3]{Furushima86}).
It is easy to see that a general hypersurface $Z_6\subset \PP(1,2,3,5)$
satisfies the above property (iv). Since $Z$ is unique up to isomorphism,
it should be isomorphic to $Z_6$.
\end{proof}

The del Pezzo surface $Z$ as above we denote by $\DP$.

\begin{stheorem}{\bf Lemma.}
\label{lemma-DuVal}
{\rm (i)} Let $Z=\DP$ 
and let $\Theta$ be a generator of 
$\Cl(Z)$. 
Then 
\[
\text{$\dim |t\Theta|=t-1$ for $t=1,\dots,4$ and $\dim |t\Theta|=t$ for $t=5$ and $6$}.
\]

{\rm (ii)} Let $Z= \PP(1,2,3)$ and let $\Theta$ be the positive generator of 
$\Cl(Z)$. 
Then 
\[
\text{$\dim |t\Theta|=t-1$ for $t=1,\dots,5$, and $\dim |6\Theta|=6$}.
\]
\end{stheorem}

\begin{proof}
Follows by the orbifold Riemann-Roch formula \cite[Theorem 9.1]{Reid-YPG1987}
and Kawamata-Viehweg vanishing.
\end{proof}

The following fact is a consequence of computer search (see \ref{computer}),
degree bound $-K_X^3\le 125/2$ \cite{Prokhorov-2007-Qe}, and 
\cite[Th. 1.4]{Prokhorov2008a}. 

\begin{mtheorem}{\bf Proposition.}
\label{proposition-dim-A}
Let $X$ be a $\QQ$-Fano threefold with $\qQ(X)=\qW(X)$ and let $A\in \Cl(X)$ be such that $-K_X=\qW(X)A$.
\begin{enumerate}
 \item 
If $\qW(X)\ge 8$, then $\dim |A|\le 0$.
 \item 
If $\qW(X)=7$ and $\dim |A|\ge 1$, then $X\simeq \PP(1^2,2,3)$.
 \item 
If $\qW(X)\ge 4$ and $X$ is singular, then $\dim |A|\le 2$.
 \item 
If $\dim |A|\ge 2$ and $\qW(X)\ge 5$, then $X\simeq \PP(1^3,2)$.
\end{enumerate}
\end{mtheorem}

\section{Case $\qQ(X)=9$}
\label{section-q=9}
Let $X$ be a $\QQ$-Fano threefold with $\qQ(X)=9$ and 
$\g(X)>4$. By \cite[Proposition 3.6]{Prokhorov2008a} we have $\B(X)=(2,4,5)$ and
$\Cl(X)\simeq \ZZ$.
Further,
$\dim |kA|=0$, $1$, $2$, $4$, $6$ for $k=1$, $2$, $3$, $4$, $5$,
respectively.
Thus there are irreducible divisors $S_k\in |kA|$ for $k=1,2,3,4,5$.

First we prove the following.
\begin{mtheorem}{\bf Proposition.}
\label{proposition-q=9-D5}
The pair $(X,|4A|)$ is canonical.
\end{mtheorem}

\begin{proof}
Put $\MMM:=|4A|$ and assume that $(X,\MMM)$ is not canonical.
Apply construction \eqref{eq3.1}.
In \eqref{equation-discrepancies} we have
$0<a_0=\beta_4-\alpha$. 
Hence, 
\[
\beta_4>\alpha,
\quad 
\beta_1\ge \frac14\beta_4> \frac14 \alpha, 
\quad
\beta_2\ge \frac12\beta_4> \frac12 \alpha.
\]
Further, for some $a_i\in \ZZ$ we can write
 \begin{equation}
\label{equation-rel-Cl-9-1}
\begin{array}{lllll}
K_{\bar X}+9\bar S_1&+ &a_1\bar E&\sim& 0,
\\[4pt]
K_{\bar X}+ \bar S_1+4 \bar S_2&+ &a_2 \bar E&\sim& 0,
\\[4pt]
K_{\bar X}+ \bar S_4+\bar S_5&+ &a_4 \bar E&\sim& 0,
\end{array}
\end{equation}
where
\begin{equation}
\label{equation-rel-Clb-1-9}
\begin{array}{lllllll}
a_1&=&9\beta_1-\alpha&>&5\beta_1&>&0,
\\ 
a_2&=&\beta_1+4\beta_2-\alpha&>&\beta_1+2\beta_2&>&0,
\\
a_4&=&\beta_4+\beta_5-\alpha&>&\beta_5&>&0.
\end{array}
\end{equation}

\begin{stheorem}
{\bf Claim.}
$f(E)$ is a point of index $>1$.
\end{stheorem}
\begin{proof}
Assume that $f(E)$ is either a Gorenstein point 
or curve.
Then $\alpha$ and $\beta_k$ are integers, so 
$a_2\ge 4$.
If $\bar f$ is not birational, then 
restricting \eqref{equation-rel-Cl-9-1} to its general fiber $V$
we get a contradiction (otherwise $-K_V$ is divisible by some integer $\ge a_2$).
Hence, the contraction $\bar f$ is birational. From the second relation of \eqref{equation-rel-Cl-9-1} 
we have $\hat q=s_1+4s_2 +a_2e\ge 8$.
Then by \cite[Th. 1.4]{Prokhorov2008a} $\hat S_2\notqq\Theta$ (because $\dim |\hat S_2|\ge 1$), so $s_2\ge 2$.
Therefore, $\hat q \ge 13$ (see \eqref{equation-index}).
Using \cite[Prop. 3.6]{Prokhorov2008a} 
we get successively $s_2\ge 3$, $\hat q\ge 17$, $s_2\ge 5$, and
$\hat q>19$, a contradiction.
\end{proof}

\begin{scase}
From now on we assume that $f(E)$ is a point of index $r>1$, 
where $r=2$, $4$, or $5$. By Lemma \ref{lemma-discrepancies} we have $\alpha=1/r$.
\end{scase}

\begin{stheorem}{\bf Claim.}
$\hat X$ is not a surface.
\end{stheorem}

\begin{proof}
Assume that $\hat X$ is a surface. 
Restricting \eqref{equation-rel-Cl-9-1} to a general fiber 
$\bar f^{-1} (\operatorname{pt})\simeq \PP^1$ we see that 
the divisors $\bar S_1$ and $\bar S_2$ are $\bar f$-vertical and $\bar S_1\sim \bar f^*\Theta$.
Since $\dim |\bar S_1|=0$, $\hat X$ is either $\PP(1,2,3)$ or $\DP$
(see \ref{case-not-birational-del-Pezzo}). 
Moreover, since $a_4>0$, $\bar S_k$ is also $\bar f$-vertical for $k=4$ or $5$.
Then $\bar S_k\sim k\bar S_1\sim k\bar f^*\Theta$.
In this case, $2k-4=\dim |\bar S_k|=\dim |k\Theta|$.
We get a contradiction by Lemma \ref{lemma-DuVal}. 
\end{proof}

\begin{stheorem}{\bf Claim.}
$\hat X$ is not a curve.
\end{stheorem}
\begin{proof}
Assume that $\hat X\simeq \PP^1$. 
Then by \eqref{equation-rel-Cl-9-1} divisors $\bar S_1$ and $\bar S_2$
are $\bar f$-vertical, $\bar S_4$ and $\bar S_5$ are $\bar f$-horizontal
(because $\bar S_4$ and $\bar S_5$ are irreducible
and $\dim |S_4|, \dim |S_5|>1$), and $\bar f$ is generically $\PP^2$-bundle.
In this case,
$3=a_1=9\beta_1-\alpha$. 
Since $\alpha=1/r$, we have 
$r\beta_1=(1+3r)/9\notin \ZZ$,
a contradiction.
\end{proof}

From now on we assume that $\bar f$ is birational.
Then
\begin{equation}
\label{equation-q=9-hatq}
9s_1+ a_1e= s_1+4s_2 +a_2e= s_4+s_5+a_4e=\hat q.
\end{equation}

\begin{stheorem}{\bf Claim.}
\label{lemma-q-9}
$\hat q\le 8$.
\end{stheorem}
\begin{proof}
Assume that $\hat q\ge 9$. Then the group $\Cl(\hat X)$ is torsion free
by Theorem \ref{theorem-main-10-19}.
Since $\alpha=1/r$, we have $\g (\hat X)\ge \g (X)=18$ by Lemma \ref{lemma-dk} and
$\hat q \le 13$ (see \cite[Table 3.6]{Prokhorov2008a}).
Again from \cite[Table 3.6]{Prokhorov2008a}, using inequalities
$\dim |\hat S_4|\ge \dim|4A|=4$ and $\dim |\hat S_5|\ge \dim |5A|=6$, one can obtain
successfully 
\begin{multline*}
\hat q\ge 9 \Longrightarrow
s_4\ge 4,\ s_5\ge 5 \Longrightarrow \hat q\ge 11
\Longrightarrow
s_4\ge 5,\ s_5\ge 6 \Longrightarrow 
\\
\Longrightarrow\hat q= 13
\Longrightarrow s_4\ge 6,\ s_5\ge 8\Longrightarrow a_4<0,
\end{multline*}
which is a contradiction.
\end{proof}

\begin{stheorem}{\bf Corollary.}
\label{corollary-q=9-sa}
$s_1=0$, $e=1$, $s_2=1$, $a_1=\hat q$, $a_2=\hat q-4$.
\end{stheorem}

\begin{stheorem}{\bf Corollary.}
$r=5$ and $\hat X\simeq \PP(1^2,2,3)$.
\end{stheorem}
\begin{proof}
By \eqref{equation-rel-Clb-1-9} we have $\hat q=a_1=9\beta_1-1/r$,
$r\beta_1=(\hat q r+1)/9$. Since $r\beta_1$ is an integer, there is only one possibility:
$r=5$ and $\hat q=7$. Hence, $s_2=1$. The group $\Cl(\hat X)$ is torsion free
by Lemma \ref{lemma-torsion-ae}.
Then the second assertion follows by Proposition \ref{proposition-dim-A}
because $\dim |\Theta|\ge 1$.
\end{proof}

Now we are in position to finish the proof of Proposition \ref{proposition-q=9-D5}.
Since $\dim |\hat S_5|\ge 6$, $s_5\ge 3$. Similarly $s_4\ge 3$.
From \eqref{equation-rel-Cl-9-1} we get $a_4=1$, 
$s_4=s_5=3$, i.e.,
$\hat S_4\sim \hat S_5\sim 3\Theta$.
Note also that $\dim |S_5|=\dim |3\Theta|=6$,
so $|S_5|=|5A|$ is the birational pull-back of $|3\Theta|$.

By Corollary \ref{corollary-q=9-sa} $\bar S_1$ is the $\bar f$-exceptional divisor and
by \eqref{equation-discrepancies-second} we have
\[
K_{\bar X}+ 7\bar S_2\equiv \bar f^*\left(K_{\hat X} +7\hat S_2\right)+ \left(b-7\gamma_2\right)\bar S_1
\ \equiv\ \left(b-7\gamma_2\right)\bar S_1.
\] 
Pushing this relation to $X$ we get $b=5+7\gamma_2\ge 5$.
In particular, $\bar f(\bar S_1)$ is a smooth point \cite{Kawamata-1996} and so
$b$, $\delta$, and $\gamma_k$ are integers (because $K_{\hat X}$, $\hat E$, and $\hat S_k$
are Cartier at $\bar f(\bar S_1)$).
Since $\Bs |3\Theta|$ is contained in the singular locus of 
$\hat X\simeq \PP(1^2,2,3)$, $\gamma_5=0$.
On the other hand,
\[
K_{\bar X}+ \frac 73\bar S_5\equiv \bar f^*\left(K_{\hat X} +\frac 73\hat S_5\right)+ \left(b-\frac 73\gamma_5\right)\bar S_1
\equiv \left(b-\frac 73\gamma_5\right)\bar S_1.
\]
Hence, $8/3=b-7\gamma_5/3=b$, a contradiction. Proposition \ref{proposition-q=9-D5} is proved completely.
\end{proof}

\begin{case}\label{explanation-DP}
Now by Proposition \ref{proposition-q=9-D5}
the pair $(X,\MMM)$ is 
canonical.
Then a general member $M\in \MMM$ is a 
normal surface with only Du Val singularities
\cite[1.12, 1.21]{Alexeev-1994ge}. By the adjunction formula
$-K_M=5A|_M$ and $K_M^2=5$.
Therefore, $M$ is a del Pezzo surface of degree $5$.
By Proposition \ref{lemma-DuVal-5} we have $M\simeq \DP\simeq M_6\subset \PP(1,2,3,5)$.
Then the statement \ref{theorem-main-q=9} of Theorem \ref{theorem-main} 
follows by the following.
\end{case}

\begin{mtheorem}{\bf Theorem (hyperplane section principle, 
\cite{Reid-graded-rings}, cf. \cite{Mori-1975}, {\cite{Sano-1996}}). }
\label{hyperplane-section-principle}
Let $R$ be a graded integral domain, and let $0\neq x_0 \in R$ be a graded element of degree $a_0$. 
Let $\bar R := R/(x_0 )$. 
\begin{enumerate}
 \item 
 Let $\bar x_1,\dots, \bar x_k$ be
homogeneous elements that generate $\bar R$, and $x_1,\dots, x_k \in R$ any 
homogeneous elements that map to $\bar x_1,\dots, \bar x_k \in\bar R$. Then $R$ is generated by
$x_0, x_1,\dots, x_k$.
 \item 
Under the assumption of \textup{(i)}, let $\bar f_ 1,\dots,\bar f_ n$ be homogeneous generators
of the ideal of relations holding between $\bar x_1,\dots, \bar x_k$. Then there exist
homogeneous relations $f_1,\dots, f_n$ holding between $x_0, x_1,\dots, x_n$ in $R$
such that the $f_i$ reduces to $\bar f_ i$ modulo $x_0$ and $f_1,\dots, f_n$ generate the
relation between $x_0, x_1,\dots, x_n$.
\end{enumerate}
\end{mtheorem}

\section{Case $\qQ(X)=6$}
Using \ref{computer} one can get the following computation.
\begin{mtheorem}{\bf Lemma.}
\label{lemma-q=6-table}
Let $X$ be a $\QQ$-Fano threefold with $\qW(X)=\qQ(X)=6$ and 
$\g(X)>15$. Let $-K_X=6A$.
Then the group $\Cl(X)$ is torsion free and we have one of the following cases:
\rm
\setlongtables\renewcommand{\arraystretch}{1.3}
\begin{longtable}{|c|c|c|c|c|c|c|c|c|}
\hline
&&&\multicolumn{6}{c|}{$\dim |kA|$}
\\
\hhline{|~|~|~|------}
&$\B$&$A^3$& $|A|$&$|2A|$&$|3A|$&$|4A|$&$|5A|$&$|-K|$
\\[5pt]
\hline
\endfirsthead
\hline
&$\B$&$A^3$& $|A|$&$|2A|$&$|3A|$&$|4A|$&$|5A|$&$|-K|$
\\[5pt]
\hline
\endhead
\hline
\endlastfoot
\hline
\endfoot

\nr\label{no-q=6-B7}&
$(7)$& $2/7$& $1$& $4$& $8$& $14$& $22$& $32$
\\
\nr\label{no-q=6-B5}&
$(5)$& $1/5$& $1$& $3$& $6$& $10$& $16$& $23$
\end{longtable}
\end{mtheorem}

\begin{case}
Let $X$ be a $\QQ$-Fano threefold with $\qQ(X)=6$ and $\g(X)>15$.
First we assume that is torsion free. General case will follow
by Lemma \ref{q=6-general-case} below.
In particular, $\qQ(X)=\qW(X)$.
Thus $X$ is described in Lemma \ref{lemma-q=6-table}.
We will show that the case \xref{no-q=6-B7}
does not occur.
The case \ref{no-q=6-B5} was considered in \cite{Sano-1996}.
It can also be treated by the method of \S \ref{section-q=9}.
\end{case}

\begin{case} \textbf{Case \xref{no-q=6-B7}.}
\label{case-q=6-b=7}
We have
\begin{equation}
\label{equation-KS-q=6}
K_{\bar X}+6\bar S_1+a_1\bar E\sim 0.
\end{equation}
Put $\MMM:=|A|$.
Near the point of index $7$ we have $A\sim -6K_X$ and $\MMM\sim A \sim -6K_X$.
By Lemma \ref{lemma-c-42} we get $c\le 1/6$.
So, $\beta_1\ge 6\alpha$ and $a_1=6\beta_1-\alpha\ge 35\alpha\ge 5$.
Therefore, the contraction $\bar f$ is birational and 
$\hat q\ge 11$. Moreover,
$f(E)\in X$ is a cyclic quotient singularity of index $7$ and 
$\alpha=1/7$ by \cite{Kawamata-1996}, Lemma \ref{lemma-discrepancies}.
Taking Lemma \ref{lemma-dk} into account we obtain $\g(\hat X)\ge \g(X)=31$.
This contradicts Theorem \ref{theorem-main-10-19}.
\end{case}

\begin{mtheorem}
\label{q=6-general-case}
{\bf Lemma.}
Let $X$ be a $\QQ$-Fano threefold with $\qQ(X)=6$.
Assume that the torsion part of $\Cl(X)$ is non-trivial.
Then $\g(X)\le 13$. 
\end{mtheorem}

\begin{proof}
Assume that $\g(X)\ge 14$. 
By Lemma \ref{lemma-torsion} there is a $\QQ$-Fano threefold $Y$
with $\Cl(Y)\simeq \ZZ$, $\g(Y)\ge 23$ and
whose Fano-Weil index $\qW(Y)$ is divisible by $6$.
This contradicts the above considered case \xref{no-q=6-B7}
 and \eqref{equation-index}. 
\end{proof}

\section{Case $\qQ(X)=8$}
\begin{mtheorem}{\bf Lemma.}
\label{lemma-q=8-table}
Let $X$ be a $\QQ$-Fano threefold with $\qQ(X)=\qW(X)=8$, $\g(X)\ge 10$ and $A^3\neq 4/91$. Let $-K_X=8A$.
Then the group $\Cl(X)$ is torsion free and we have one of the following cases:
\rm
\setlongtables\renewcommand{\arraystretch}{1.3}
\begin{longtable}{|c|c|c|c|c|c|c|c|c|c|c|}
\hline
&&&\multicolumn{8}{c|}{$\dim |kA|$}
\\
\hhline{|~|~|~|--------}
&$\B$&$A^3$& $|A|$&$|2A|$&$|3A|$&$|4A|$&$|5A|$&$|6A|$&$|7A|$&$|-K|$
\\[5pt]
\hline
\endfirsthead
\hline
&$\B$&$A^3$& $|A|$&$|2A|$&$|3A|$&$|4A|$&$|5A|$&$|6A|$&$|7A|$&$|-K|$
\\[5pt]
\hline
\endhead
\hline
\endlastfoot
\hline
\endfoot

\nr\label{no-q=8-3-9}&
$(3, 9)$ &
$1/9
$ &$0
$ &$2
$ &$4
$ &$7
$ &$11
$ &$16
$ &$22
$ &$29$

\\
\nr\label{no-q=8-3-5-11-a}&
$(3, 5, 11)$ &
$16/165
$ &$0
$ &$1
$ &$3
$ &$5
$ &$9
$ &$13
$ &$18
$ &$25$

\\
\nr\label{no-q=8-11}&
$(11)$ &
$1/11
$ &$0
$ &$2
$ &$3
$ &$6
$ &$9
$ &$13
$ &$18
$ &$24$

\\
\nr\label{no-q=8-7-11}&
$(7, 11)$ &
$6/77
$ &$0
$ &$1
$ &$2
$ &$4
$ &$7
$ &$10
$ &$15
$ &$20$

\\
\hline
&&&&&&&&&&

\\[-12pt]

\nr\label{no-q=8-3-3-5}&
$(3, 3, 5)$ &
$1/15
$ &$0
$ &$1
$ &$3
$ &$4
$ &$7
$ &$10
$ &$13
$ &$18$

\\

\nr\label{no-q=8-3-7}&
$(3, 7)$ &
$1/21
$ &$0
$ &$1
$ &$2
$ &$3
$ &$5
$ &$7
$ &$10
$ &$13$

\\
\hline
&&&&&&&&&&

\\[-12pt]

\nr\label{no-q=8-3-3-5-9}&
$(3, 3, 5, 9)$ &
$2/45
$ &$-1
$ &$0
$ &$1
$ &$2
$ &$4
$ &$6
$ &$8
$ &$11$

\\
\hline
%
%
\end{longtable}
\end{mtheorem}

\begin{case}
Now let $X$ be a $\QQ$-Fano threefold with $\qQ(X)=8$ and $\g(X)>10$.
First we assume that the group $\Cl(X)$ is torsion free. For general case see
Lemma \ref{q=8-general-case} below.
Thus $X$ is described in Lemma \ref{lemma-q=8-table}.
We will show that cases
\xref{no-q=8-3-9}, 
\xref{no-q=8-3-5-11-a},
\xref{no-q=8-11}, 
\xref{no-q=8-7-11}, and \xref{no-q=8-3-3-5-9}
do not occur. We also show that in the cases
 \xref{no-q=8-3-3-5} and \xref{no-q=8-3-7}
we have $X\simeq X_6\subset \PP(1,2,3^2,5)$ and $X\simeq X_{10}\subset \PP(1,2,3,5,7)$,
respectively.
In all cases, $\dim |A|\le 0$ and linear systems 
$|2A|$ and $|3A|$ have no fixed components.
Write
\begin{equation}
\label{equation-KS-q=8}
\begin{array}{lll}
K_{\bar X}&+8\bar S_1&+a_1\bar E\sim 0,
\\[5pt]
K_{\bar X}&+4\bar S_2&+a_2\bar E\sim 0,
\\[5pt]
K_{\bar X}&+2\bar S_3+\bar S_2&+a_3\bar E\sim 0,
\end{array}
\end{equation}
where $a_i\in \ZZ$.
 \end{case}

\begin{case}{\bf Cases \xref{no-q=8-3-9} and \xref{no-q=8-11}.}
\label{q=8-small-table-13}
Put $\MMM:=|2A|$ and let $P\in X$ be a (unique) point
of index $9$ (resp. $11$) in the case \xref{no-q=8-3-9} (resp. \xref{no-q=8-11}).
Then $\MMM\sim -mK_X$ near $P$, where $m=7$ (resp. $m=3$) 
in the case \xref{no-q=8-3-9} (resp. \xref{no-q=8-11}).
Hence, $c\le 1/m$ (see \ref{lemma-c-42}).
Thus $\beta_2\ge m\alpha$ and $\beta_1\ge \beta_2/2\ge m\alpha/2$.
In particular, 
\[
a_1\ge 8\beta_1-\alpha\ge (4m-1)\alpha,\quad a_2=4\beta_2-\alpha\ge (4m-1)\alpha.
\]

\begin{stheorem}
{\bf Claim.} 
The contraction $\bar f$ is birational.
\end{stheorem}
\begin{proof}
Assume that $\bar f$ is a contraction of fiber type.
The divisors $\bar S_2$ and $\bar S_1$ are $\bar f$-vertical by 
\eqref{equation-KS-q=8}.
In particular, $\bar S_2\sim 2\bar S_1$.
Since $\dim |\bar S_2|=2$ and $\bar S_2$ is irreducible,
$\hat X$ is not a curve. By \ref{case-not-birational-del-Pezzo}
we have $\hat X\simeq \PP(1^2,2)$ or $\DP$.
This contradicts Lemma \ref{lemma-DuVal}.
\end{proof}

\begin{stheorem}
{\bf Claim.} 
 $f(E)$ is a point of index $>1$.
\end{stheorem}
\begin{proof}
Assume that $f(E)$ is either a curve or Gorenstein point.
Then all the numbers $\alpha$ and $\beta_l$ are integers.
We have 
$a_2\ge 11$ and $\hat q\ge 4s_2+a_2\ge 15$. In this case, $|\Theta|=\emptyset$
and $\dim |2\Theta|\le 0$ (see Theorem \ref{theorem-main-10-19}).
Therefore, $s_2\ge 3$ and $\hat q>19$, a contradiction.
\end{proof}

\begin{scase}
Thus we may assume that the contraction $\bar f$ is birational and $f(E)$ is a point of index $r>1$.
One can see from Table \ref{lemma-q=8-table} that $f(E)$ is a cyclic quotient singularity.
In particular, $\alpha=1/r$ \cite{Kawamata-1996}.
Since $\alpha<1$, we have $\g(\hat X)\ge \g (X)=23$ by Lemma \ref{lemma-dk}.
By Theorem \ref{theorem-main-10-19} and \ref{theorem-main-q=9} of Theorem \ref{theorem-main}
we also have $\hat q\le 8$.
By \eqref{equation-KS-q=8} we can write
\begin{equation}
\label{equation-KS-q=8-1}
\begin{array}{lll}
\hat q=8s_1+a_1e=4s_2+a_2e=2s_3+s_2+a_3e.
\end{array}
\end{equation}
Since $a_1,\, a_2\ge 1$, we have $s_1=0$, $s_2=1$, and $\hat q\ge 5$. 
By Lemma \ref{lemma-torsion-ae} the group $\Cl(\hat X)$ is torsion free. 
By Proposition \xref{proposition-dim-A} 
$\hat X\simeq \PP(1^3,2)$.
Therefore, $1=a_2=4\beta_2-1/r=3\beta_2$, $\beta_2=1/3$, $r=3$, and we are 
in the case \xref{no-q=8-3-9}. 
This contradicts $\beta_2\ge m\alpha=7/3$.
\end{scase}
\end{case}

\begin{case} {\bf Case \xref{no-q=8-3-5-11-a}.}
\label{q=8-small-table}
Put $\MMM:=|3A|$ 
and let $P\in X$ be a (unique) point of index $11$.
Then $\MMM\sim -10K_X$ near $P$. 
Then $c\le 1/10$ (see \ref{lemma-c-42}).
Thus $\beta_3\ge 10\alpha$ and $\beta_1\ge 10\alpha/3$.
Further,
\[
a_1\ge 77\alpha/3,\quad 
a_3\ge 19\alpha+\beta_2
\quad \Longrightarrow \quad a_1\ge 3,\quad a_3\ge 2.
\]

\begin{stheorem}
{\bf Claim.} 
The contraction $\bar f$ is birational.
\end{stheorem}
\begin{proof}
Assume that $\bar f$ is a contraction of fiber type.
Since $a_1\ge 3$, $\hat X$ is a curve and $\bar f$ is a generically $\PP^2$-bundle.
On the other hand, $\bar S_3$ is $\bar f$-vertical by 
\eqref{equation-KS-q=8} and because $a_3\ge 2$. 
Since $\dim |\bar S_3|=3$ and $\bar S_3$ is irreducible, 
we get a contradiction.
\end{proof}

\begin{stheorem}
{\bf Claim.} 
 $f(E)$ is a point of index $>1$.
\end{stheorem}
\begin{proof}
Assume that $f(E)$ is either a curve or Gorenstein point.
Then all the numbers $\alpha$ and $\beta_l$ are integers.
We have 
$a_3\ge 19\alpha\ge 19$. Thus, $\hat q\ge 2s_3+a_3>19$, a contradiction.
\end{proof}
\begin{scase}
Thus we may assume that the contraction $\bar f$ is birational and $f(E)\in X$ is a 
a cyclic quotient singularity of index $r=3$, $5$
or $11$ (see Table \ref{lemma-q=8-table}).
In particular, $\alpha=1/r$ \cite{Kawamata-1996}. Hence, $\g(\hat X)\ge \g (X)=23$ by Lemma \ref{lemma-dk}.
Consider the equality \eqref{equation-KS-q=8-1}.
By Theorem \ref{theorem-main-10-19} and \ref{theorem-main-q=9} of Theorem \ref{theorem-main}
we have $\hat q\le 8$.
Note that $\beta_2\ge 1/r$ (because $S_2$ cannot be Cartier at $f(E)$).
Hence, $a_2=4\beta_2-1/r>0$.
By \eqref {equation-KS-q=8-1} we have $s_1=0$, $s_2=1$, and $\hat q\ge 5$. 
By Lemma \ref{lemma-torsion-ae} the group $\Cl(\hat X)$ is torsion free and $e=1$. 
Since $s_2=1$, $\dim |\Theta|\ge \dim |2A|\ge 1$. 
By Proposition \ref{proposition-dim-A} we have $s_3\ge 2$, so $\hat q\ge 7$.
Again applying Proposition \ref{proposition-dim-A} we get 
$\hat X\simeq \PP(1^2,2,3)$.
Hence, $7=\hat q=a_1=8\beta_1-1/r$ and so
$r\beta_1=(7r+1)/8$. One can see that this number is not integral for 
$r=3$, $5$ and $11$,
a contradiction. 
\end{scase}
\end{case}

\begin{case} \textbf{Case \xref{no-q=8-7-11}.}
Take $\MMM:=|5A|$. Near the point of index $7$ we have $\MMM\sim -5K_X$.
Then in Lemma \ref{lemma-c-42} we have $m=5$, so 
$\beta_5\ge 5\alpha$ and $\beta_1\ge \alpha$. Thus $a_1\ge 7\alpha$.

\begin{stheorem}
{\bf Claim.} 
 $f(E)$ is a point of index $r>1$.
\end{stheorem}
\begin{proof}
Assume that $f(E)$ is either a curve or Gorenstein point.
Then $\alpha \in \ZZ$, $a_1\ge 7$. Therefore, the contraction $\bar f$ is birational and 
$\hat q\ge 8s_1+7e$. If $s_1>0$, then $\hat q\ge 17$, $|\Theta|=\emptyset$,
$e\ge 2$, and $\hat q>19$, a contradiction. Hence, $s_1=0$.
Then by Lemma \ref{lemma-torsion-ae} the group $\Cl(\hat X)$ is torsion free and $e=1$.
Moreover, $\hat q=a_1=8\beta_1-\alpha=8\alpha+7(\beta_1-\alpha)$. 
Since the last term is non-negative, by \eqref{equation-index} we get $\hat q=7$ and $\alpha=1$.
Similar to \eqref{equation-KS-q=8} write 
\begin{equation}
\label{equation-q=8-expr-dob}
K_{\bar X}+\bar S_3+\bar S_5+a_5\bar E\sim 0,\quad a_5=\beta_3+\beta_5-\alpha.
\end{equation}
Thus $a_5\ge 4$ and $s_3+s_5\le \hat q-a_5\le 3$. 
By Proposition \ref{proposition-dim-A} we have $\hat X\simeq \PP(1^2,2,3)$, so $\dim |3\Theta|=6$.
This contradicts $\dim |5A|=7$.
\end{proof}

\begin{scase}\label{subcase-q=8-first-m}
\textbf{Subcase: $r=7$.}
Then $\alpha=1/7$ and $A\sim -K_X$ near $f(E)$. 
Hence, for $l=1,\dots,6$, we have $\beta_l=l/7+m_l$, where $m_l$ is a non-negative integer.
This gives us $a_1=8\beta_1-\alpha= 1+8m_1$.
If $\bar f$ is not birational, then 
restricting \eqref{equation-KS-q=8} to a general fiber $V$ one can see 
that $\bar S_1$ is $\bar f$-vertical, $a_1=1$, and $-K_V=\bar E|_{V}$. 
In particular, $\hat X\not \simeq \PP^2$, $\PP(1^2,2)$. 
By \eqref{equation-q=8-expr-dob} we have $a_5=1+m_3+m_5$, so
$a_5=1$ and the divisors $\bar S_3$ and $\bar S_5$ are $\bar f$-vertical.
This contradicts Lemma \ref{lemma-DuVal}.

Hence the contraction $\bar f$ is birational. Using \eqref{equation-KS-q=8} we write
\[
\hat q=e+ 8(s_1+m_1e)= e+ 4(s_2+m_2e)
\]
Since $s_2>0$, we have $\hat q\ge 9$. 
Since $\alpha<1$, we have $\g(\hat X)\ge \g (X)=19$ by Lemma \ref{lemma-dk}.
From Theorem \ref{theorem-main-10-19} and \ref{theorem-main-q=9} of Theorem \ref{theorem-main}
we conclude that $\hat q=11$ and $\hat X\simeq \PP(1,2,3,5)$.
Further, $e=3$ and \eqref{equation-q=8-expr-dob} gives us
\[
11=\hat q= s_3+s_5+a_5e=3+ s_3+s_5+3m_3+3m_5,\quad s_3+s_5\le 8.
\]
On the other hand, $\dim |\hat S_3|\ge 2$ and $\dim |\hat S_5|\ge 7$, so 
$s_3\ge 4$ and $s_5\ge 6$, a contradiction.
\end{scase} 

\begin{scase}
\textbf{Subcase: $r=11$.}
Then $\alpha=1/11$ and $A\sim -7K_X$ near $f(E)$. 
Thus $\beta_1=7/11+m_1$ and $\beta_5=2/11+m_5$, where the $m_i$ are 
non-negative integers. Since $\beta_5\ge 5\alpha$, $m_5\ge 1$.
Further,
$a_1=5+8m_1$, $a_5=1+m_1+m_5$. Hence, the contraction $\bar f$ is birational and
\[
\hat q=5e+8(s_1+m_1e)= e+ s_3+s_5+m_3e+m_5e. 
\]
By \eqref{equation-index} either $\hat q=13$ or $\hat q=5$.
Since $\g(\hat X)\ge \g (X)=19$, we have $\hat q=5$.
Thus $s_1=0$, $e=1$, $\Cl(\hat X)$ is torsion free, and $s_3+s_5\le 3$.
Therefore, $s_3=1$.
By 
Proposition \ref{proposition-dim-A} $\hat X\simeq \PP(1^3,2)$.
In this case, $\dim |2\Theta|=6$. This contradicts 
$s_5\le 2$ because $\dim |5A|=7$.
\end{scase}
\end{case}

\begin{case} \textbf{Case \xref{no-q=8-3-3-5-9}.}
 Put $\MMM:=|3A|$.
Near the point of index $9$ we have $A\sim -8K_X$, $\MMM\sim -6K_X$.
Thus $c\le 1/6$ and $\beta_3\ge 6\alpha$.
We can write
\begin{equation}
\label{equation-KS-q=8-7-13}
\begin{array}{lll}
K_{\bar X}&+4\bar S_2&+a_2\bar E\sim 0,
\\[5pt]
K_{\bar X}&+2\bar S_3+S_2&+a_3\bar E\sim 0. 
\end{array}
\end{equation}
\[
a_3= 2\beta_3+\beta_2-\alpha\ge 11\alpha+\beta_2.
\]
If $\alpha\ge 1$, then $a_3\ge 11$. Hence the contraction $\bar f$ is birational.
In this case, $\hat q\ge 11e+2s_3\ge 13$,
$s_3\ge 2$, $\hat q\ge17$, $e>1$, and $q>19$, a contradiction.
Therefore, $\alpha<1$ and $P:=f(E)$ is a point of index $r>1$.

\begin{scase} \textbf{Subcase $r=5$.} 
Then near $P$ we have $-K_X\sim 3A$ and $A\sim -2K_X$.
Hence, $\beta_2=4/5+m_2$ and $\beta_3=1/5+m_3$, where $m_3\ge 1$.
We get 
\[
a_2=3+4m_2,\qquad a_3=1+2m_3+m_2. 
\]
If $\bar f$ is not birational, then it is a generically
$\PP^2$-bundle and $a_2=a_3=3$. In this case, $\bar S_3$ 
is a general fiber. On the other hand, $\bar S_2$ is $\bar f$-vertical.
Hence, $m\bar S_2$ is also a scheme fiber for some $m\in \ZZ_{>0}$.
This implies $\bar S_3\sim m\bar S_2$ and $S_3\sim mS_2$,
a contradiction.
Therefore, the contraction $\bar f$ is birational. Then
\[
\hat q=4(s_2+m_2e)+3e=2(s_3+m_3e) +s_2+m_2e+e\ge 5,
\]
\[
3(s_2+m_2e)=2(s_3+m_3e-e).
\]
Write $s_2+m_2e=2k$, $s_3+m_3e-e=3k$, where $k\in \ZZ_{>0}$.
Then $\hat q=8k+3e\ge 11$ and $\Cl(\hat X)$ is torsion free.
By Lemma \ref{lemma-torsion-ae} $F\sim eA$, so $e>1$,
$\hat q\ge 17$, $s_3\ge 5$, $k\ge 2$, and $q>19$,
a contradiction.
\end{scase}

\begin{scase} \textbf{Subcase $r=3$ and $\alpha=1/3$.} 
Then near $P$ we have $-K_X\sim 2A$ and $A\sim -2K_X$.
Hence, $\beta_2=1/3+m_2$ and $\beta_3=m_3$, where $m_3\ge 2$.
We get 
\[
a_2=1+4m_2,\qquad a_3=2m_3+m_2\ge 4. 
\]
Therefore, the contraction $\bar f$ is birational. In this situation,
\[
19\ge \hat q=4(s_2+m_2e)+e=2(s_3+m_3e) +s_2+m_2e\ge 6,
\]
Hence, $s_2+m_2e\le 4$.
If $s_2+m_2e$, then so $\hat q$ is.
In this case, $\hat q\le 8$ and so $s_2+m_2e=0$.
Then $e=\hat q\ge 6$ and 
$\hat q=2(s_3+m_3e)>19$, a contradiction.
If $s_2+m_2e=1$, then $\hat q$ is odd, $e=\hat q-4\ge 3$,
$\hat q=2(s_3+m_3e)+1\ge 15$, $s_3\ge 5$, and $\hat q>19$. Again we get a contradiction.

Finally, assume that $s_2+m_2e=3$. Then 
$\hat q\ge 13$ and $\Cl(\hat X)$ is torsion free.
Since $|A|=\emptyset$, we have $>1$ 
by Lemma \ref{lemma-torsion-ae}. Hence, $\hat q\ge 17$
and so $s_3\ge 5$, $\hat q>19$, a contradiction.
\end{scase}

\begin{scase} \textbf{Subcase $r=3$ and $\alpha=2/3$.} 
Then near $P$ we have $-K_X\sim 2A$ and $A\sim -2K_X$.
Hence, $\beta_2=2/3+m_2$ and $\beta_3=m_3$, where $m_3\ge 4$.
We get 
\[
a_2=2+4m_2,\qquad a_3=2m_3+m_2\ge 8. 
\]
Therefore, the contraction $\bar f$ is birational. Then
\[
19\ge \hat q=4(s_2+m_2e)+2e=2(s_3+m_3e) +s_2+m_2e\ge 6.
\]
Hence $\hat q $ is even and $\hat q\ge 10$. This contradicts \eqref{equation-index}.
\end{scase}

\begin{scase} \textbf{Subcase $r=9$.} 
Then $\beta_2=7/9+m_2$, $\beta_3=6/9+m_3$,
$a_2=3+4m_2$, $a_3=2+2m_3+m_2$. If $\bar f$ is not birational,
then it is a generically $\PP^2$-bundle, $a_2=3$,
$m_2=0$, $m_3=0$, $a_3=2$, $\bar S_3$ is a general fiber, and
$\bar S_2$ is $\bar f$-vertical.  Hence, 
$\OOO_{\bar S_3}(a_2\bar E)=\OOO_{\bar S_3}(a_3\bar E)=\OOO_{\bar S_3}(-K_{\bar X})$,
 a contradiction.

Therefore, $\bar f$ is not birational and so
\[
\hat q=4(s_2+m_2e)+3e=2(s_3+m_3e) +s_2+m_2e+2e.
\]
If $s_2+m_2e\ge 2$, then $\hat q\ge 11$, the group $\Cl(\hat X)$ is torsion free,
and $e\ge 2$ by Lemma \ref{lemma-torsion-ae}. In this case, 
both $s_2+m_2e$ and $e$ are odd $\ge 3$. Hence, $\hat q>19$,
 a contradiction. Therefore, $s_2+m_2e\le 1$.

If $s_2+m_2e=0$, then 
\[
\hat q=3e=2(s_3+m_3e)+2e, \quad \hat q=6,
\quad e=2, \quad s_3=1, \quad F=S_2.
\]
By Lemma \ref{lemma-torsion-ae} the group $\Cl(\hat X)$ is torsion free.
Using computer search with 
$\dim |\Theta|\ge 1$ and 
\ref{theorem-main-q=6} of Theorem \ref{theorem-main} we get
$X\simeq X_6\subset \PP(1^2,2,3,5)$.
Thus we have $K_{\bar X}+4\bar S_2+3\bar E\sim 0$
and $K_{\bar X}+2\bar S_3+\bar S_2+2\bar E\sim 0$.
Further, $K_{\bar X}+2\bar S_4+a_4\bar E\sim 0$,
where $a_4=2\beta_4-\alpha$, $\beta_4=5/9+m_4$, and
$a_4=1+2m_4$. 
We get $\hat q=6=2s_4+2em_4+e=2s_4+4m_4+2$,
$s_4=2$.
Comparing dimensions of linear systems we see that
$|\Theta|$ is the birational transform of $|3A|$
and $|2\Theta|$ is the birational transform of $|4A|$.
Write
\[
\begin{array}{l}
K_{\bar X}+6\bar S_3=\bar f^*(K_{\hat X}+6\hat S_3)+(b-6\gamma_3)\bar F,
\\[5pt]
K_{\bar X}+3\bar S_4=\bar f^*(K_{\hat X}+3\hat S_4)+(b-3\gamma_4)\bar F.
\end{array}
\]
Since $F\sim 2 A$, we have $b=5+6\gamma_3\ge 5$ and 
$b=2+3\gamma_4$.
Hence, $Q:=\bar f(\bar F)$ is a Gorenstein point.
Assume that $\gamma_3=0$. 
Let $\hat S_3'\in |\Theta|$ be a (unique)
divisor passing through $Q$.
Then $\hat S_3' \neq \hat E$ (because $e>1$) and
\[
K_{\bar X}+6\bar S_3'=\bar f^*(K_{\hat X}+6\hat S_3')+(b-6\gamma_3')\bar F. 
\]
where $\gamma_3'>0$. Let $S_3'\sim mA$.
Then $m>2$ (because $|A|=\emptyset$ and 
$S_3'\neq S_2=F$). So, $-8+6m=2(b- 6\gamma_3')$ and $m=3- 2\gamma_3'\le 1$, a contradiction.
We get $\gamma_3>0$ and $\gamma_4=1+2\gamma_3>0$.
Therefore, $Q\in \Bs |\Theta|\cap \Bs |2\Theta|$ and so $Q$ is the point of index
$5$, a contradiction.

If $s_2+m_2e=1$, then $\hat q =4+3e=2(s_3+m_3e)+2e+1$.
If furthermore $e\ge 2$, then $s_2=1$, $\hat q\neq 19$,
$\hat q =13$, $e=3$, $s_3+m_3e=3$, and $s_3=3$. Moreover, the linear system
$|3\Theta|$ is the birational transform of $|3A|$.
On the other hand, by Lemma \ref{lemma-torsion-ae} 
the map $X\dashrightarrow \hat X$ contracts a divisor $F\sim 3A$,
a contradiction.
Therefore, $\hat q=7$, $e= 1$, and $s_3+m_3=2$.
Since $|A|=\emptyset$, by Lemma \ref{lemma-torsion-ae} we have 
$\Cl(\hat X)\simeq \ZZ\oplus \ZZ_n$, where $n>1$.
If $s_2=1$, then $\hat E\qq \hat S_2$.
In this situation,  by \cite[Th. 1.4]{Prokhorov2008a} $\hat X\simeq \PP(1^2,2,3)$
and $\Cl(\hat X)\simeq \ZZ$, a contradiction.
Therefore, $s_2=0$, i.e., $S_2=F$.
By Lemma \ref{lemma-q=7-torsion-general} we have 
$\qQ(\hat X)=\qW(\hat X)$. Since $\dim |-K_{\hat X}|\ge 11$ and 
$\Cl(\hat X)\simeq \ZZ\oplus \ZZ_n$, 
using computer search and \ref{theorem-main-q=7} of Theorem \ref{theorem-main} we get 
only one possibility: $\B(\hat X)=(2, 2, 3, 12)$, $A^3=1/12$.
In this case, $n=2$ or $4$ \cite[Prop. 2.9]{Prokhorov2008a}, 
so there is an \'etale in codimension two 
cover $\pi : X'\to \hat X$ of degree $n$ such that $K_X'=\pi^*K_X$.
Moreover, $X'$ is a Fano threefold with terminal singularities,
$\qW(X')=\qQ(X')=7$, and $\B(X')=(3,3,6)$ or $(3,3,3,3,3)$.
In particular, the Gorenstein index of $X'$ is strictly smaller than
its Fano index. According to \cite{Sano-1996} there is only one 
such a Fano threefold whose Gorenstein index equals to $6$:
$\PP(1,1,2,3)$. Then $X'$ does not have any point of index $6$, a contradiction.
For $\B(X')=(3,3,3,3,3)$ we have two possibilities:
$X'_6\subset \PP(1^2,2,3^2)$ and $X'_4\subset \PP(1^3,2,3)$ in both cases 
$\qW(X')=4$. Again this is a contradiction. 
\end{scase}
\end{case}

\begin{case} \textbf{Case \xref{no-q=8-3-3-5}.}
We show that in this case $X$ is a hypersurface of degree $6$ in $\PP(1,2,3^2,5)$.
Similar to \ref{explanation-DP} it is sufficient to prove the following.

\begin{stheorem}{\bf Lemma.}
Let $X$ be a $\QQ$-Fano threefold with $\qW(X)=8 $
and $\B(X)=(3,3,5)$. Then the pair $(X,|3A|)$ is canonical. 
\end{stheorem}
\begin{proof}
Put $\MMM:=|3A|$ and assume that $(X,\MMM)$ is not canonical. 
Then $\beta_3>\alpha$, $\beta_1\ge \frac13 \beta_3>\frac13 \alpha$, and so
\begin{equation}
\label{equation-q=8-A4}
a_1=8\beta_1-\alpha> \frac 53 \alpha>0, \qquad
a_3=2\beta_3+\beta_2-\alpha>\alpha+\beta_2>0.
\end{equation}

\begin{stheorem}
{\bf Claim.} 
 $f(E)$ is a point of index $r>1$.
\end{stheorem}
\begin{proof}
Assume that $f(E)$ is either a curve or Gorenstein point.
Then as in \eqref{equation-q=8-A4} we have 
$a_1\ge 5\beta_1+\beta_3-\alpha\ge 5\beta_1+1\ge 6$ and
$a_3\ge \beta_3+\beta_2+1\ge 3$.
Therefore, 
 the contraction $\bar f$ is birational and 
\[
\hat q=2s_3+s_2+a_3e= 8s_1+a_1e\ge 6.
\]
If $\hat q\ge 9$, then $\Cl(\hat X)\simeq \ZZ$ by Theorem \ref{theorem-main-10-19} and  
using \cite[Table 3.6]{Prokhorov2008a}
we obtain successively $s_3\ge 4$, $\hat q\ge 13$, $s_3\ge 5$, $\hat q\ge 17$, 
$s_3>7$, and $\hat q>19$, a contradiction. Therefore, 
$\hat q \le 8$. In this case, $s_1=0$, and $\Cl(\hat X)\simeq \ZZ$ by Lemma \ref{lemma-torsion-ae}.
By Proposition \ref{proposition-dim-A} we have
$\dim |\Theta|<3$, so $s_3>1$ and $\hat q\ge 8$.
Therefore, $\hat q= 8$, $s_3\le 2$,
and $s_2= 1$. This contradicts Proposition \ref{proposition-dim-A}.
Thus $f(E)$ is a point of index $r=3$ or $5$.
Then $\alpha<1$ by Lemma \ref{lemma-discrepancies}. 
Since $\g(\hat X)\ge \g(X)=17$, we have $\hat q\le 13$.
\end{proof}

\begin{scase}
{\bf Subcase $r=5$.}
Then $\alpha=1/5$ and $A\sim -2K_X$ near $f(E)$.
So,
\[
\beta_1=2/5+m_1,\quad \beta_2=4/5+m_2,\quad \beta_3=1/5+m_3, 
\]
where $m_i$ are non-negative integers and $m_3\ge 1$ (cf. \ref{subcase-q=8-first-m}). 
Hence,
\[
a_1=3+8m_1,\quad a_2=3+4m_2,\quad a_3=1+m_2+2m_3.
\]
If $\bar f$ is not birational, then $\hat X\simeq \PP^1$
and the divisors $\bar S_1$, $\bar S_2$, $\bar S_3$ are $\bar f$-vertical.
This contradicts $\dim | \bar S_3|=3$.
Therefore, the contraction $\bar f$ is birational. 
We have
\[
\hat q=8s_1+3e+8m_1e=4s_2+ 3e+4m_2e=2s_3+s_2+e+m_2e+2m_3e.
\]
Hence either $s_1=0$ or $\hat q\ge 11$.
In both cases $\Cl(\hat X)\simeq \ZZ$ (see 
Theorem \ref{theorem-main-10-19} and Lemma \ref{lemma-torsion-ae}).
Further, $\hat q \ge 7$ and $s_3\ge 2$ by Proposition \ref{proposition-dim-A}. 
If $e>1$, then $\hat q \ge 11$, $s_3\ge 4$, $\hat q\ge 17$, 
$s_3>7$, $\hat q>19$, a 
contradiction. Hence, $e=1$ and $\hat q\equiv 3\mod 8$.
This is possible only if $\hat q=11$. Then
\[
s_1+m_1=1, \quad
s_2+m_2=2,\quad
s_3+m_3=4,\quad s_3\le3.
\]
This contradicts $\dim |3\Theta|\le 2$.
\end{scase}
\begin{scase}
{\bf Subcase $r=3$.} 
First, we assume that $\alpha=1/3$.
Then $A\sim -2K_X$ near $f(E)$.
So,
\[
\beta_1=2/3+m_1,\quad \beta_2=1/3+m_2,\quad \beta_3=m_3, 
\]
where $m_i$ are non-negative integers and $m_3\ge 1$ (cf. \ref{subcase-q=8-first-m}). 
Hence,
\[
a_1=5+8m_1,\quad a_2=1+4m_2,\quad a_3=m_2+2m_3.
\]
Since $a_1>3$, the contraction $\bar f$ is birational. We have
\[
\hat q=8s_1+5e+8m_1e=4s_2+ e+4m_2e=2s_3+s_2+m_2e+2m_3e.
\]
We have either $s_1=0$ or $\hat q\ge 13$.
In both cases $\Cl(\hat X)\simeq \ZZ$ (see Lemma \ref{lemma-torsion-ae}).
Further, $\hat q \ge 5$ and $s_3\ge 2$ by Proposition \ref{proposition-dim-A}.
Hence, $\hat q\ge 7$.
If $e>1$, then $\hat q \ge 11$, $s_2\ge 2$, $s_3\ge 4$, and $\hat q\ge 14$, a 
contradiction. 
Hence, $e=1$ and $\hat q\equiv 5\mod 8$.
This is possible only if $\hat q=13$. Then
\[
s_1+m_1=1, \quad
s_2+m_2=3,\quad
s_3+m_3=5,\quad s_3\le4.
\]
This contradicts $\dim |4\Theta|\le 2$.

Finally we consider the case where 
$\alpha\ge 2/3$. Then
$P:=f(E)$ is a point of index $3$ with basket $(3,3)$.
By Lemma \ref{lemma-discrepancies} we have
$\alpha=2/3$.
Further, $A\sim -2K_X$ and $3A\sim 0$ near $f(E)$.
So, $\beta_1=4/3+m_1$, where $m_1$ is a non-negative integer, and 
$a_1=10+8m_1>3$.
Hence, the morphism $\bar f$ is birational. In this case,
\[
13\ge \hat q=10e+8s_1+8m_1e.
\]
The only possibility is $\hat q=10$.
This contradicts \eqref{equation-index}.
\end{scase}
\end{proof}
\end{case}

\begin{case} \textbf{Case \xref{no-q=8-3-7}.}
We will show that in this case $X\simeq X_{10}\subset\PP(1,2,3,5,7)$.
Put $\MMM:=|2A|$.
Near the point of index $7$ we have $A\sim -K_X$ and $\MMM\sim -2K_X$.
Thus $c\le 1/2$, $\beta_2\ge 2\alpha$, $\beta_1\ge \alpha$, and 
\begin{equation}
\label{equation-KS-q=8-2}
a_1=8\beta_1-\alpha\ge 7\alpha,\quad
a_2=4\beta_2-\alpha\ge 7\alpha.
\end{equation}
If $\alpha\ge 1$, then $a_1,\, a_2\ge 7$, $\bar f$ is birational, 
and $\hat q\ge 4s_2+7\alpha\ge 11$.
Hence, $s_2\ge 2$,
$\hat q\ge 17$, $s_2\ge 5$, and $\hat q>19$, a contradiction.
Thus $\alpha<1$. Therefore, 
$f(E)$ is a point of index $r=3$ or $7$, and $\alpha=1/r$.

\begin{scase} \textbf{Subcase $r=3$.}
Then $\beta_2=1/3+m_2$,
where $m_2\ge 1$.
By \eqref{equation-KS-q=8-2} we have
$a_2=1+4m_2$.
Hence the contraction $\bar f$ is birational.
Then
\[
\hat q=e+4(s_2+m_2e).
\]
Since $s_2>0$, we have $\hat q\ge 9$
and so $\Cl(\hat X)\simeq \ZZ$. We get successively $s_2\ge 2$,
$\hat q\ge 13$, $s_2\ge 3$, $\hat q\ge 17$, $s_2\ge 5$, and $\hat q>19$,
a contradiction.
\end{scase}

\begin{scase} \textbf{Subcase $r=7$.}
Then, as above, $\beta_k=k/7+m_k$ for $k=1,\dots,6$, so
$a_1=1+8m_1$.
If $\bar f$ is not birational, then $\bar S_1$ 
is $\bar f$-vertical. 
Restricting \eqref{equation-KS-q=8} to a general fiber $V$ we get 
$-K_V\sim \bar E|_{V}$ and $a_1=1$.
Then using relations similar to \eqref {equation-KS-q=8}
we get that in this situation $\bar S_k$ 
is $\bar f$-vertical for $k=1,\dots,6$. 
Since $\dim |\bar S_3|=2$, $\hat X$ is not a curve.
Since $\dim |\bar S_1|=0$, $\hat X$ is either $\PP(1,2,3)$ of $\DP$. 
On the other hand, $\dim |\bar S_6|= 7$.
This contradicts Lemma \ref{lemma-DuVal}.
Hence the contraction $\bar f$ is birational and then
\[
\hat q=e+8(s_1+m_1e)=e+4(s_2+m_2e),\quad 2(s_1+m_1e)=s_2+m_2e\ge s_2>0.
\]
Hence, $\hat q\ge 9$ and $s_2\ge 2$.
If $\hat q\ge 13$, then $s_2\ge 3$, $\hat q\ge 17$, $s_2\ge 5$, and $\hat q>19$, a contradiction.
Thus $\hat q\le 11$ and so 
\[
s_1+m_1e=1,\quad
\hat q=e+8,\quad s_2+m_2e=2,\quad s_2=2,\quad m_2=0.
\]

Assume that $s_1=0$, i.e., $\bar F=\bar S_1$. 
Then $e=m_1=1$ by Lemma \ref{lemma-torsion-ae}. So, $\hat q=9$,
and $\hat X\simeq X_6\subset \PP(1,2,3,4,5)$
by Theorem \ref{theorem-main}, \ref{theorem-main-q=9}.
Recall that $c=\alpha/\beta_2=1/2$. By our construction \eqref{eq3.1} the pair $(\hat X,\frac12 \hat \MMM)$ 
is canonical, where $\hat \MMM\subset |2\Theta|$.
On the other hand, near
the point $\hat P_5\in \hat X$ of index $5$ we have 
$\hat \MMM\sim -3K_{\hat X}$. 
By Lemma \ref{lemma-c-42} we have 
$\ct(\hat X,\hat \MMM)\le 1/3$, a contradiction.

Therefore, $s_1=1$. Then $e>1$. So $\hat q= 11$, $e=3$, 
and $\hat X\simeq \PP(1,2,3,5)$ by \cite{Prokhorov2008a}.
Similar to \eqref{equation-KS-q=8} we get $s_k=k$ for $k=1,\dots,6$.
Since $\dim |k\Theta|=\dim |kA|$,
the linear system $|k\Theta|$ is the birational transform of $|kA|$ for $k=1,\dots,6$.
Write
\begin{equation}
\label{equation-q=8-711-Kgamma}
K_{\bar X}+\frac{11}{k}\bar S_k\equiv
\bar f^*\left(K_{\hat X}+\frac{11}{k}\hat S_k\right)
+\left(b- \frac{11}{k}\gamma_k\right)\bar F,
\end{equation}
where $K_{\hat X}+\frac{11}{k}\hat S_k\equiv 0$. 
Note that $F\sim 3A$ by Lemma \ref{lemma-torsion-ae}.
Pushing \eqref{equation-q=8-711-Kgamma} down to $X$ we get 
\[
3= 3\left(b- \frac{11\gamma_k}{k}\right), \quad 
b= 1+ \frac{11\gamma_k}{k}\ge 1.
\]
Since $\hat X=\PP(1,2,3,5)$ has only cyclic quotient singularities, 
$\bar f(\bar F)$ is either a curve or a smooth point
\cite{Kawamata-1996}.
On the other hand, $|6\Theta|$ is base point free outside of the point of index $5$.
Therefore, $\gamma_6=0$, $b=1$, and $C:=\bar f(\bar F)$ is a curve.
Recall that $\hat E\in |\OOO_{\PP(1,2,3,5)}(3)|$ and so $\hat E\simeq \PP(1,2,5)$.
Since $C$ is 
contained into the smooth locus of $\hat X$ \cite{Kawamata-1996},
the curve $C$ is cut out on $\hat E$ by a 
divisor $\hat S\in |\OOO_{\PP(1,2,3,5)}(l)|$, where $l$ is divisible by $10$. 
On the other hand, $\bar E\simeq \hat E$ is covered by a family of curves 
$L_\lambda\in |\OOO_{\PP(1,2,5)}(2)|$ and
\[
0\le -K_{\bar X}\cdot L_\lambda = -K_{\hat X}\cdot 
\bar f(L_\lambda)- \bar F \cdot L_\lambda=\frac{11}5-\frac l 5. 
\]
Hence, $l=10$.

Further, the restriction $-K_{\bar X}|_{\bar E} =\OOO_{\PP(1,2,5)}(1)$
is ample. Note that
\begin{equation*}
-K_{\bar X}\equiv \frac{11}{3}\bar E+\frac{8}{3}\bar F
\equiv
\frac{11}{6}\bar S_6
-\bar F.
\end{equation*}
Since the only base point of $|\hat S_6|=|6\Theta|$
is the point of index $5$, the divisor $-K_{\bar X}$ is ample.
Therefore, the map $\bar X\dashrightarrow \tilde X$ is a composition of 
flips with respect to $K$ (or an isomorphism).
On the other hand, $\bar X$ and $\tilde X$ have the same collections 
of non-Gorenstein singularities: $\frac12(1,1,1)$, $\frac 13 (1,1,2)$,
and $\frac 15 (1,2,3)$. 
This means that the map $\bar X\dashrightarrow \tilde X$ does not change 
Shokurov's difficulty \cite{Shokurov-1985}. Hence 
$\bar X\dashrightarrow \tilde X$ is an isomorphism.
Since $\hat S_1\cdot C=\hat S_1\cdot\hat E\cdot\hat S_{10}=1$, the intersection 
$\hat S_1\cap C$ is a single smooth point, say $\hat P$.
The diagram \eqref{eq3.1} induces the following 
\begin{equation*}
\begin{gathered}
\xymatrix{
&\bar S_1\ar[dr]^{\bar h}\ar[dl]_{h}&
\\
S_1&&\hat S_1
} 
\end{gathered}
\end{equation*}
where $\hat S_1\simeq \PP(2,3,5)$, $\bar h$ is the blowup of
$\hat P$, and $h$ is the contraction of the proper transform of 
the curve $\hat S_1\cap \hat E$ given by $x_3=0$.
This shows that the surface $S_1$ is unique up to isomorphism.
On the other hand, taking $X= X_{10}\subset\PP(1,2,3,5,7)$ we get 
the same $S_1$.
Therefore, for any choice of $X$ of type \xref{no-q=8-3-7}
the surface $S_1\in |A|$ is isomorphic to 
a general weighted hypersurface of degree $10$ in $\PP(2,3,5,7)$.
Then by Theorem \ref{hyperplane-section-principle} we conclude that 
$X\simeq X_{10}\subset\PP(1,2,3,5,7)$.
\end{scase}
\end{case}

\begin{mtheorem}
\label{q=8-general-case}
{\bf Lemma.}
Let $X$ be a $\QQ$-Fano threefold with $\qQ(X)=8$.
Assume that the torsion part of $\Cl(X)$ is non-trivial.
Then $\g(X)\le 8$. 
\end{mtheorem}

\begin{proof}
Similar to Lemma \ref{q=6-general-case}.
\end{proof}

\section{Case $\qQ(X)=7$}
\begin{mtheorem}{\bf Lemma.}
 \label{lemma-q=7-torsion-general}
Let $X$ be a $\QQ$-Fano threefold with $\qQ(X)=7$. Then 
$\qW(X)=7$.
\end{mtheorem}
\begin{proof}
Assume that $\qQ(X)\neq \qW(X)$. Then there is a 
$7$-torsion element $\Xi\in \Cl(X)$ \cite[Lemma 3.2 (iv)]{Prokhorov2008a}.
Let $\B^{\Xi}\subset \B(X)$
is the subbasket of points where $\Xi$ is not Cartier.
According to \cite[Prop. 2.9]{Prokhorov2008a} $\B^{\Xi}=(7,7,7)$.
Kawamata's condition $\sum_{P\in \B(X)} (r_P-1/r_P)<24$ gives us that 
$\B(X)\setminus \B^{\Xi}=\emptyset$, $(2)$, $(2,2)$, or $(3)$. 
The element $\Xi$ defines 
a cyclic \'etale in codimension two cover $\pi: X'\to X$ of degree $7$ such that $X'$ is a Fano threefold with 
terminal singularities and $\qQ(X')$ is divisible by $7$. 
On the other hand, $X'$ is Gorenstein at points $\pi^{-1}(\B^{\Xi})$.
Hence the Gorenstein index of $X'$ is at most $3$.
This contradicts \cite{Sano-1996}.
\end{proof}

\begin{mtheorem}{\bf Lemma.}
\label{lemma-q=7-table}
Let $X$ be a $\QQ$-Fano threefold with $\qQ(X)=7$ and let $-K_X=7A$.
Assume either $\g(X)> 17$ or $\g(X)= 17$ and $A^3=1/10$.
Then the group $\Cl(X)$ is torsion free and we have one of the following cases:
\rm
\setlongtables\renewcommand{\arraystretch}{1.3}
\begin{longtable}{|c|c|c|c|c|c|c|c|c|c|c|}
\hline
&&&\multicolumn{8}{c|}{$\dim |kA|$}
\\
\hhline{|~|~|~|--------}
&$\B$&$A^3$& $|A|$&$|2A|$&$|3A|$&$|4A|$&$|5A|$&$|6A|$&$|-K|$
\\[5pt]
\hline
\endfirsthead
\hline
&$\B$&$A^3$& $|A|$&$|2A|$&$|3A|$&$|4A|$&$|5A|$&$|6A|$&$|-K|$
\\[5pt]
\hline
\endhead
\hline
\endlastfoot
\hline
\endfoot

\nr\label{no-q=7-B23}&
$(2, 3)$
&$1/6
$&$1
$&$3
$&$6
$&$10
$&$15
$&$22
$&$30$

\\
\hline
&&&&&&&&&

\\[-12pt]

\nr\label{no-q=7-B2259}&
$(2, 2, 5, 9)$
&$7/45
$&$0
$&$2
$&$4
$&$8
$&$13
$&$19
$&$27$

\\
\nr\label{no-q=7-B228}&
$(2, 2, 8)$
&$1/8
$&$0
$&$2
$&$4
$&$7
$&$11
$&$16
$&$22$

\\
\nr\label{no-q=7-B2368}&
$(2, 3, 6, 8)$
&$1/8
$&$0
$&$1
$&$3
$&$6
$&$10
$&$16
$&$22$

\\
\nr\label{no-q=7-B2359}&
$(2, 3, 5, 9)$
&$11/90
$&$0
$&$1
$&$3
$&$6
$&$10
$&$15
$&$21$

\\
\nr\label{no-q=7-B23410}&
$(2, 3, 4, 10)$
&$7/60
$&$0
$&$1
$&$3
$&$6
$&$9
$&$14
$&$20$

\\
\hline
&&&&&&&&&

\\[-12pt]
\nr\label{no-q=7-B2225}&
$(2, 2, 2, 5)$
&$1/10
$&$0
$&$2
$&$3
$&$6
$&$9
$&$13
$&$18$

\end{longtable}
\end{mtheorem}

\begin{case}
Now let $X$ be a $\QQ$-Fano threefold with $\qQ(X)=7$.
Assume either $\g(X)>20$ or $\g(X)=17$ and $A^3=1/10$.
By Lemma \ref{lemma-q=7-torsion-general}
$\qQ(X)=\qW(X)$. 
Thus $X$ is described in Lemma \ref{lemma-q=7-table} and $\Cl(X)$ is torsion free.
The case \ref{no-q=7-B23} was considered in \cite{Sano-1996} and \cite{Prokhorov2008a}.
We will show that cases
\ref{no-q=7-B2259} -- \ref{no-q=7-B23410}
do not occur and in the case 
\ref{no-q=7-B2225}
we have $X\simeq X_6\subset \PP(1,2^2,3,5)$.
Write 
\begin{equation}
\label{equation-q=7-B=2-3-6-8}
\begin{array}{lllll}
K_{\bar X}&+7\bar S_1&&+a_1\bar E&\sim 0,
\\[5pt]
K_{\bar X}&+\bar S_1 &+3\bar S_2 &+a_2\bar E&\sim 0.
\end{array}
\end{equation}
where $a_1$, $a_2$ are some integers.
\end{case}

\begin{case}
First we consider cases \ref{no-q=7-B2259} -- \ref{no-q=7-B23410}.
Put $\MMM:=|2A|$.
Let $P\in X$ be a point of index $r'>0$, where we take $r'$ as follows:
\begin{equation}
\label{equation-q=7-rm}
\begin{array}{c|ccccc}
\text{No.}
&\text{\ref{no-q=7-B2259}}
&\text{\ref{no-q=7-B228}}
&\text{\ref{no-q=7-B2368}}
&\text{\ref{no-q=7-B2359}}
&\text{\ref{no-q=7-B23410}}
\\
r' & 9 & 8 & 8&9&10
\\ 
m & 8 & 6 & 6&8&6
\end{array}
\end{equation}
Here $m$ is an integer such that $0\le m<r'$
and $\MMM\sim -mK_X$ near $P$. 
Then $c\le 1/m$ (see \ref{lemma-c-42}), so 
$\beta_2\ge m\alpha$ and $\beta_1\ge \beta_2/2\ge m\alpha/2$.
We have
\begin{equation}
\label{equation-q=7-B=2-3-6-8-a}
\begin{array}{l}
a_1=7\beta_1-\alpha,
\\[5pt]
a_2=\beta_1+3\beta_2-\alpha.
\end{array}
\end{equation}
In particular, 
\begin{equation}
\label{equation-q=7-B=2-3-6-8-ab}
a_1,a_2\ge (7 m/2-1)\alpha\ge 20\alpha\ge 2. 
\end{equation}

\begin{stheorem}
{\bf Claim.} 
We have $\alpha<1/3$. Moreover,
$f(E)$ is a cyclic quotient singularity of index $r>3$ and $\alpha=1/r$.
\end{stheorem}
\begin{proof}
If $\alpha\ge 1/3$, then $a_2\ge 7$. 
In this case, the contraction $\bar f$ is birational and 
$\hat q\ge 3s_2+7$. So, 
$\hat q\ge 11$.
Since $\dim |\Theta|\le 0$, we get successively $s_2\ge 2$, $\hat q\ge 13$,
$s_2\ge 3$, $\hat q\ge 17$, $s_2\ge 5$, and $\hat q>19$, a contradiction.
Therefore, $\alpha< 1/3$. From Table \ref{lemma-q=7-table}
we infer that $f(E)\in X$ is a
cyclic quotient singularity of index $r>3$ and $\alpha=1/r$ \cite{Kawamata-1996}.
\end{proof}

\begin{stheorem} {\bf Claim.}\label{claim-q=7-non-birational}
If the morphism $\bar f$ is not birational, then 
it is a generically $\PP^2$-bundle. Moreover, $\bar S_2$ is a general fiber, 
$\bar S_2\sim 2\bar S_1$, and 
$\OOO_{\bar S_2}(\bar E)=\OOO_{\bar S_2}(\bar S_3)=\OOO_{\PP^2}(1)$. 
This is possible only in cases \xref{no-q=7-B2368} -- \xref{no-q=7-B23410}.
\end{stheorem}
\begin{proof}
Assume that $\bar f$ is not birational. Since
$a_1,\, a_2>0$, the divisors $\bar S_1$ and $\bar S_2$ 
are $\bar f$-vertical. Since $\dim |\bar S_1|=0$,
$\hat X$ is either $\PP^1$, $\PP(1,2,3)$ or $\DP$.
By Lemma \ref{lemma-DuVal} the divisor $\bar S_2$ 
is $\bar f$-horizontal.
Similar to \eqref{equation-q=7-B=2-3-6-8}
write
\begin{equation*}
K_{\bar X}+\bar S_1 +2\bar S_3 +a_3\bar E\sim 0,
\qquad a_3=\beta_1+2\beta_3-\alpha>0.
\end{equation*}
Hence, $a_3=1$ and $\bar f$ is a generically $\PP^2$-bundle.
Since $\bar S_2$ is irreducible, $\dim |\bar S_2|\le 1$.
The rest is obvious. 
\end{proof}

\begin{stheorem} {\bf Claim.}
\label{claim-q7-}
Assume that $X$ is of type \xref{no-q=7-B2259}, \xref{no-q=7-B228} or \xref{no-q=7-B2368}.
Then the morphism $\bar f$ is birational,
$X$ is of type \xref{no-q=7-B2368},
and $r=6$ or $8$.
Moreover, 
$\hat X\simeq \PP(1,2,3,5)$, $\PP(1^2,2,3)$ or 
$X_6\subset \PP(1^2,2,3,5)$.
\end{stheorem}

\begin{proof}
If $\bar f$ is not birational, then
we are in the case \ref{no-q=7-B2368} and
$a_1,\, a_2\le 3$. Hence, $a_1=a_2=3$
by \eqref{equation-q=7-B=2-3-6-8-ab}. 
By \eqref{equation-q=7-B=2-3-6-8-a} we have $\alpha\le 3/20$, so $r=8$ and $\alpha=1/8$.
Then $\beta_1=m\alpha+m_1=7/8+m_1$,
where $m_1$ is a non-negative integer (cf. \ref{subcase-q=8-first-m}).
Again by \eqref{equation-q=7-B=2-3-6-8-a} we obtain $a_1=6+7m_1\ge 6$, a contradiction.

Therefore, the contraction $\bar f$ is birational. In this case,
\[
\hat q=7s_1+a_1e=s_1 +3s_2 +a_2e\ge 6.
\]
By Lemma \ref{lemma-dk} $\g(\hat X)\ge \g(X)\ge 21$.
Hence, by Theorem \ref{theorem-main-10-19} 
and \ref{theorem-main-q=9}-\ref{theorem-main-q=8} of Theorem \ref{theorem-main},
we have either $6\le \hat q\le 7$ or $\hat q=11$ (and $\hat X\simeq \PP(1,2,3,5)$).
Moreover, $\Cl(\hat X)$ is torsion free by Lemma \ref{lemma-torsion}.
If $\hat X\simeq \PP(1,2,3,5)$, then
$s_2\le 2$ and so $\dim |2A|\le |2\Theta|=1$. 
Hence, we are in the case \ref{no-q=7-B2368}. 

Assume that $6\le \hat q\le 7$. Then $s_2=1$.
If $\hat q=7$ (resp. $\hat q=6$), then $\hat X\simeq \PP(1^2,2,3)$
(resp. $\hat X\simeq X_6\subset \PP(1^2,2,3,5)$)
by Proposition \ref{proposition-dim-A} (resp. by 
\ref{theorem-main-q=6} of Theorem \ref{theorem-main}).
Hence, $\dim |2A|\le \dim |\Theta|=1$ and we are again in the case \ref{no-q=7-B2368}.
\end{proof}

\begin{scase} {\bf Subcase \ref{no-q=7-B2368} with $r=6$.}
Then
$\beta_1=1/6+m_1$, and $\beta_2=2/6+m_2$, 
where $m_i$ are positive integers.
We get
\[ 
\hat q\ge a_1=7\beta_1-\alpha =1+7m_1\ge 8,\quad 
a_2=1+\beta_1+3\beta_2-\alpha=2+m_1+3m_2.
\]
Hence, $\hat q=11$, $\hat X\simeq \PP(1,2,3,5)$ by Claim \ref{claim-q7-}, and 
$s_2\ge 2$. Therefore,
\[
11=\hat q=2e+ s_1+m_1e+3(s_2+m_2e)\ge 2+m_1+3(s_2+m_2)\ge 12,
\]
a contradiction.
\end{scase}

\begin{scase} {\bf Subcase \ref{no-q=7-B2368} with $r=8$.}
Then 
$\beta_1=7/8+m_1$ and $\beta_2=6/8+m_2$,
where $m_i$ are non-negative integers.
We obtain
\[ 
a_1=6+7m_1\ge 6,\quad 
a_2=3+m_1+3m_2.
\]
So,
\[
\hat q=6e+7(s_1+m_1e)=3e+ s_1+m_1e+3(s_2+m_2e).
\]
Since $\hat q\in \{6, 7, 11\}$, we have $s_1=m_1=0$, $\hat q=6$, and so $\hat X\simeq 
X_6\subset \PP(1^2,2,3,5)$ by Claim \ref{claim-q7-}.
In particular, 
$\dim |\Theta|=1$,
$\dim |2\Theta|=3$,
$\dim |3\Theta|=6$, and $\Cl(\hat X)\simeq \ZZ$. Furthermore,
$e=1$ and $s_2=1$.

By our computations $\bar S_1$ is the $\bar f$-exceptional divisor,
$a_1=6$, and $a_2=3$.
By \eqref{equation-discrepancies-second} 
\[
K_{\bar X}+6\bar S_2\sim \bar f^*\left( K_{\hat X}+6\hat S_2 \right)+ \left( b-6\gamma_2\right) \bar S_1\equiv
\left( b-6\gamma_2\right) \bar S_1.
\]
Therefore, 
\[
b=6\gamma_2+5\ge 5.
\]
In particular, $\bar f(\bar S_1)$ is a Gorenstein point.
Similar to \eqref{equation-q=7-B=2-3-6-8} we can write
\begin{equation}\label{equation-q=7-B=2-3-6-8a}
\begin{array}{lllll}
K_{\bar X}&+\bar S_3&+4\bar S_1 &+a_3\bar E&\sim 0,
\\[5pt]
K_{\bar X}&+\bar S_4&+3\bar S_1 &+a_4\bar E&\sim 0.
\end{array}
\end{equation}
Easy computation as above give us $a_3\ge 4$ and $a_4\ge 3$.
Hence, $s_3\le 2$ and $s_4\le 3$.
On the other hand, comparing Table \ref{lemma-q=7-table}, 
\ref{no-q=7-B2368} and Table \ref{lemma-q=6-table}, 
\ref{no-q=6-B5} we see $s_3= 2$, $a_3= 4$, $s_4= 3$, and $a_4=3$.
Moreover, $\dim |3A|=\dim |2\Theta|$ and $\dim |4A|=\dim |3\Theta|$.
Therefore, $|2\Theta|$ (resp. $|3\Theta|$) is the birational transform 
of $|3A|$ (resp. $|4A|$).
Again using \eqref{equation-discrepancies-second} 
 we get
\[
5\le b=3\gamma_3+2=2\gamma_4+1\quad \Longrightarrow \quad \gamma_3,\, \gamma_4>0.
\]
This means that 
$\bar f(\bar S_1)$ is contained into $\Bs |2\Theta|$ and $\Bs |3\Theta|$.
The only point satisfying these conditions is the point of index $5$,
a contradiction.
\end{scase}

\begin{scase} {\bf Subcase \ref{no-q=7-B2359} with $r=5$.}
Then $\beta_1=3/5+m_1$ and $\beta_2=1/5+m_2$,
where $m_i$ are non-negative integers and $m_2\ge 2$. 
Further, 
\[ 
a_2=\beta_1+3\beta_2-\alpha=1+m_1+3m_2\ge 7.
\]
This implies that $\bar f$ is birational
and 
\[
\hat q=e+ s_1+m_1e+3(s_2+m_2e)\ge 10.
\]
Hence, $\hat q\ge 11$. 
In this case, $\Cl(\hat X)\simeq \ZZ$ and $\dim |\Theta|\le 0$, so 
$s_2\ge 2$, $\hat q \ge 13$, $s_2\ge 3$, $\hat q \ge 17$, 
$s_2\ge 5$, and $\hat q >19$,
a contradiction.
\end{scase}

\begin{scase} {\bf Subcase \ref{no-q=7-B2359} with $r=9$.}
Then $\beta_1=4/9+m_1$, and $\beta_2=8/9+m_2$,
where $m_i$ are non-negative integers. 
Further, 
\[ 
a_1=7\beta_1-\alpha =3+7m_1,\quad 
a_2=\beta_1+3\beta_2-\alpha=3+m_1+3m_2\ge 3.
\]
If $\bar f$ is birational, then
\[
\hat q=3e+7(s_1+m_1e)=3e+ s_1+m_1e+3(s_2+m_2e).
\]
Then $2(s_1+m_1e)=s_2+m_2e>0$.
Since $\hat q\neq 10$, we have 
$\hat q\ge 13$.
In this case, $19\ge \g(\hat X)\ge \g(X)=20$ by Lemma \ref{lemma-dk},
a contradiction.

Finally assume that $\bar f$ is not birational. Then
we are in the situation of Claim \ref{claim-q=7-non-birational}.
Similar to \eqref{equation-q=7-B=2-3-6-8}
write
\[ 
K_{\bar X}+\bar S_4+3\bar S_1+a_4\bar E\sim 0, \quad a_4=\beta_4+3\beta_1-\alpha=
2+m_4,\quad m_4\in \ZZ_{\ge 0}.
\]
Hence, $a_4=2$, $\bar S_4$ is $\bar f$-horizontal, and $\bar S_4\sim \bar E+4\bar S_1$.
Moreover, $\OOO_{\bar S_1}(\bar S_3)=\OOO_{\bar S_1}(\bar E)=\OOO_{\PP^2}(1)$. 
From the exact sequence 
\[
0 \longrightarrow 
\OOO_{\bar X}(\bar S_4-\bar S_2)
\longrightarrow
\OOO_{\bar X}(\bar S_4)
\longrightarrow
\OOO_{\bar S_2}(\bar S_4)
\longrightarrow 0.
\]
we get $\dim H^0(\OOO_{\bar X}(\bar S_4-\bar S_2))\ge 
7-3=4$. Therefore, $\dim |2A|=\dim |S_4-S_2|\ge 3$, a contradiction.
\end{scase}

\begin{scase} {\bf Subcase \ref{no-q=7-B23410} with $r=4$.}
Then $\beta_1=3/4+m_1$, and $\beta_2=1/2+m_2$,
where $m_i$ are non-negative integers and $m_2\ge 1$. 
Further, 
\[ 
a_1=7\beta_1-\alpha =5+7m_1,\quad a_2=\beta_1+3\beta_2-\alpha=2+m_1+3m_2\ge 5.
\]
This immediately implies that $\bar f$ is birational
and 
\[
\hat q=5e+7(s_1+m_1e)= 2e+ s_1+m_1e+3(s_2+m_2e)\ge 8.
\]
If $s_1+m_1e>0$, then $\hat q \ge 17$, $|\Theta|=\emptyset$,
$e>1$, $\hat q = 17$, $e=2$, $s_1+2m_1=1$, and $s_1=1$, a contradiction.
Hence, $s_1=m_1=0$ and $\hat q=5e$. This contradicts 
Theorem \ref{theorem-main-10-19} \ref{theorem-main-q=10}.
\end{scase}

\begin{scase} {\bf Subcase \ref{no-q=7-B23410} with $r=10$.}
Then $\beta_1=3/10+m_1$, $\beta_2=6/10+m_2$, and $\beta_3=9/10+m_2$,
where $m_i$ are non-negative integers. 
Further, in \eqref{equation-q=7-B=2-3-6-8} and \eqref{equation-q=7-B=2-3-6-8a} we have
\[ 
a_1=2+7m_1,\quad a_2=2+m_1+3m_2,\quad
a_3=2+m_1+2m_3,\quad m_i\in \ZZ_{\ge 0}.
\]
Assume that $\bar f$ is not birational. Then $a_1=a_2=a_3=2$, and $\bar S_1$, $\bar S_2$, $\bar S_3$ are 
$\bar f$-vertical. Hence $\hat X$ is not a curve. 
This contradicts Lemma \ref{lemma-DuVal}.
Therefore, the contraction $\bar f$ is birational. Then
\[
\hat q-2e=7(s_1+m_1e)= s_1+m_1e+3(s_2+m_2e)= s_1+m_1e+2(s_3+m_3e) \ge 3.
\]
If $s_1+m_1e=0$, then $s_2+m_2e=0$. This is impossible because 
$\dim |\hat S_2|>0$. Since $\hat q\le 19$ and $\hat q\neq 16$, we have 
$s_1+m_1e=1$ and $\hat q=7+2e\ge 9$.
In this case, $\Cl(\hat X)\simeq \ZZ$. Further, 
$s_3+m_3e=3$. On the other hand, since $\dim |3\Theta|\le 2$,
we have $s_3\ge 4$, a contradiction.

\end{scase}
\end{case}

\begin{case}
\textbf{Case \ref{no-q=7-B2225}.}
We show that in this case $X$ is a hypersurface of degree $6$ in $\PP(1,2^2,3,5)$.
Similar to \ref{explanation-DP} it is sufficient to prove the following. 

\begin{stheorem}{\bf Lemma.}
Let $X$ be a Fano threefold with $\qW(X)=7$ and $\B(X)=(2,2,2,5)$.
Then the pair $(X,|2A|)$ is canonical.
\end{stheorem}

\begin{proof}
Put $\MMM:=|2A|$ and 
assume that $(X,\MMM)$ is not canonical. 
In \eqref{equation-discrepancies} we have $0<a_0=\beta_2-\alpha$. 
Hence, $\beta_2>\alpha$ and $\beta_1\ge \frac12 \beta_2> \frac12 \alpha$.
Further, 
\begin{equation}
\label{equation-rel-Clb-1-79-prop}
a_1=7\beta_1-\alpha>5\beta_1>0,
\quad
a_2=\beta_1+3\beta_2-\alpha>\beta_1+2\beta_2>0.
\end{equation}
If $\bar f$ is not birational, then $\bar S_1$ and $\bar S_2$ are $\bar f$-vertical.
Since $\dim |\bar S_2|=2$, $\hat X$ is not a curve.
Since $\dim |\bar S_1|=0$, $\hat X$ is either $\PP(1,2,3)$ or $\DP$.
We have a contradiction by Lemma \ref{lemma-DuVal} because $\dim |\bar S_2|=2$.
Therefore, the contraction $\bar f$ is birational. We can write
\[
\hat q=7s_1+a_1e=s_1+3s_2+a_2e.
\]
Assume that $\beta_2\ge 2$. Then $\beta_1\ge \frac12 \beta_2\ge 1$, $a_2\ge 6$, and $\hat q\ge 9$.
In this case, $\Cl(\hat X)\simeq \ZZ$,
$\dim |2\Theta|\le 1$, $s_2\ge 3$, $\hat q\ge 17$, $s_2\ge 7$, and $\hat q\ge 28$,
a contradiction. Therefore, $\alpha<\beta_2<2$.
If $f(E)$ is either a Gorenstein point or curve, then $\alpha, \, \beta_2\in \ZZ$,
and so $\alpha\le 0$. Again we get a contradiction.
Thus $f(E)$ is a point of index $r=2$ or $5$.

\begin{scase}
{\bf Subcase $r=2$.} Then $\beta_2\in \ZZ$, so $\beta_2=1$ and $\alpha=1/2$.
Since $A\sim -K_X$ near $f(E)$, we have $\beta_1=1/2+m_1$,
where $m_1$ is a non-negative integer.
As above we have 
\[
a_1=3+7m_1,\quad a_2=m_1+3,
\]
\[
\hat q=7(s_1+m_1e)+3=s_1+m_1e+3(s_2+e)\ge 6.
\]
The only possibility is $s_1+m_1e=2$, $\hat q=17$, $s_2+e=5$. 
But then $\dim |6\Theta|=1$, so $s_2\ge 7$, a contradiction.
\end{scase}
\begin{scase}
{\bf Subcase $r=5$.}
Then $\alpha=1/5$, $\beta_1=3/5+m_1$, and $\beta_2=1/5+m_2$, where
$m_1$ and $m_2$ are non-negative integers with $m_2\ge 1$.
Since $\g(\hat X)\ge \g(X)=17$,
we have $\hat q\le 13$.
Further,
\[
a_1=4+7m_1,\quad a_2=1+m_1+3m_2\ge 4.
\]
This gives us
\[
\hat q = 4e+ 7(s_1+m_1e)=e+s_1+m_1e+3(s_2+m_2e)\ge 7.
\]
Hence, $\hat q=8$ or $11$. If $\hat q=8$, then $s_1=m_1=0$
and $e=2$. On the other hand, $e=1$ by Lemma \ref{lemma-torsion-ae}, a contradiction.
Therefore, $\hat q=11$, $e=1$, $s_1+m_1=1$, $s_2+m_2=3$, and $s_2\le 2$.
On the other hand, 
$\dim |2\Theta|<\dim |2A|$, a contradiction.
\end{scase}
\end{proof}
\end{case}

\section{Case $\qW(X)=5$}

\begin{mtheorem}{\bf Lemma.}
\label{lemma-q=5-table}
Let $X$ be a $\QQ$-Fano threefold with $\qW(X)=\qQ(X)=5$. Let
$-K_X=5A$.
Assume that $\g(X)\ge 19$.
Then the group $\Cl(X)$ is torsion free and we have one of the following cases:
\rm
\setlongtables\renewcommand{\arraystretch}{1.3}
\begin{longtable}{|c|c|c|c|c|c|c|c|}
\hline
&&&\multicolumn{5}{c|}{$\dim |kA|$}
\\
\hhline{|~|~|~|-----}
&$\B$&$A^3$& $|A|$&$|2A|$&$|3A|$&$|4A|$&$|-K|$
\\[5pt]
\hline
\endfirsthead
\hline
&$\B$&$A^3$& $|A|$&$|2A|$&$|3A|$&$|4A|$&$|-K|$
\\[5pt]
\hline
\endhead
\hline
\endlastfoot
\hline
\endfoot
\nr\label{no-q=5-wps}&
$(2)$&
$1/2$& $2 $& $6 $& $12 $& $21 $& $ 33$

\\
\hline
&&&&&&&

\\[-12pt]

\nr\label{no-q=5-B2236}&
$(2, 2, 3, 6)$&
$1/2$& $1$& $5$& $11$& $20$& $32$

\\
\nr\label{no-q=5-B67}&
$(6, 7)$&
$19/42$& $1$& $4$& $10$& $18$& $29$

\\ \nr\label{no-q=5-B246}&
$(2, 4, 6)$&
$5/12$& $1$& $4$& $9$& $17$& $27$

\\ \nr\label{no-q=5-B237}&
$(2, 3, 7)$&
$17/42$& $1$& $4$& $9$& $16$& $26$

\\ \nr\label{no-q=5-B228}&
$(2, 2, 8)$&
$3/8$& $1$& $4$& $8$& $15$& $24$

\\ \nr\label{no-q=5-B79}&
$(7, 9)$&
$22/63$& $1$& $3$& $7$& $13$& $22$

\\
\hline
&&&&&&&

\\[-12pt]
\nr\label{no-q=5-B223}&
$(2, 2, 3)$&
$1/3$&$1$&$4$&$8$&$14$&$22$

\\
\hline
&&&&&&&

\\[-12pt]
\nr\label{no-q=5-B47}&
$(4, 7)$&
$9/28$&$1$&$3$&$7$&$13$&$21$

\\ \nr\label{no-q=5-B412}&
$(4, 12)$&
$1/3$&$1$&$3$&$7$&$13$&$21$

\\ \nr\label{no-q=5-B222336}&
 $(2, 2, 2, 3, 3, 6)$&
$1/3$&$0$&$3$&$7$&$13$&$21$


\end{longtable}
\end{mtheorem}

\begin{case}
Let $X$ be a $\QQ$-Fano threefold with $\qQ(X)=5$.
Assume that $\g(X)\ge 21$.
We also assume that $\qQ(X)=\qW(X)$. The case $\qQ(X)\neq\qW(X)$ is excluded by 
Lemma \ref{q=5-general-case} below. 
Thus $X$ is described in Lemma \ref{lemma-q=5-table} and $\Cl(X)$ is torsion free
by Lemma \ref{lemma-torsion}.
The case \ref{no-q=5-wps} is considered in \cite{Sano-1996}.
We will show that cases \xref{no-q=5-B2236}
-- \xref{no-q=5-B79} and \xref{no-q=5-B47} -- \xref{no-q=5-B222336}
do not occur and in the case
 \xref{no-q=5-B223} we have $X\simeq X_4\subset \PP(1^2,2^2,3)$.
Note that the case \xref{no-q=5-B2236} was disproved in \cite{Prokhorov-2007-Qe}.

In cases \xref{no-q=5-B2236} -- \xref{no-q=5-B79} we put $\MMM:=|A|$.
Let $P\in X$ be a point of index $r'$, where we take $r'$ as follows:
\[
\begin{array}{c|cccccccccc}
\text{No.}
&\text{\xref{no-q=5-B2236}}
&\text{\xref{no-q=5-B67}}
&\text{\xref{no-q=5-B246}}
&\text{\xref{no-q=5-B237}}
&\text{\xref{no-q=5-B228}}
&\text{\xref{no-q=5-B79}}
&\text{\xref{no-q=5-B47}}
&\text{\xref{no-q=5-B412}}
\\
r' & 6& 6& 6& 7&8&7& 7 & 12 &
\\ 
m & 5& 5& 5& 3&5&3 & 3 & 5 &
\end{array}
\]
Here $m$ is an integer such that $0\le m<r'$
and $\MMM\sim -mK_X$ near $P$. 
Then $c\le 1/m$ (see \ref{lemma-c-42}) and so $\beta_1\ge m\alpha$.

For some $a_i\in \ZZ$ we can write
 \begin{equation}
\label{equation-rel-Cl-79-1}
-K_{\bar X} \sim 5\bar S_1+a_1\bar E \sim \bar S_1+2\bar S_2+a_2\bar E,
\end{equation}
where
\begin{equation}
\label{equation-rel-Clb-1-79}
a_1=5\beta_1-\alpha\ge (5m-1)\alpha,
\quad
a_2=\beta_1+2\beta_2-\alpha.
\end{equation}
\end{case}

\begin{case} \textbf{Cases \ref{no-q=5-B2236}, \ref{no-q=5-B67}, \ref{no-q=5-B246},
\ref{no-q=5-B228}.}
By \eqref{equation-rel-Clb-1-79} we have $a_1\ge 3$.
Assume that the contraction $\bar f$ is birational.
Since $\hat q> a_1\ge 24\alpha$, we have $\alpha<1$.
Hence, $\g(\hat X)\ge \g (X)\ge 23$.
On the other hand, $\hat q\ge 5s_1+a_1\ge 9$.
This contradicts 
Theorems
\ref{theorem-main-10-19}
and
\ref{theorem-main}
\ref{theorem-main-q=9}-\ref{theorem-main-q=8}.

Thus $\bar f$ is not birational. Then $3\ge a_1\ge 24\alpha$. 
This is possible only in the case \ref{no-q=5-B228} and then $a_1=3$,
$\bar f$ is a 
generically $\PP^2$-bundle and $\bar S_1$ is $\bar f$-vertical.
Thus $\bar S_1\simeq \PP^2$ is a general fiber.
Furthermore, $\alpha=1/8$. 
This means that $f(E)$ is the point of index $8$ and
$\beta_1=5/8$. 
Similar to \eqref{equation-rel-Cl-79-1} write
\[
-K_{\bar X} \sim \bar S_3+2\bar S_1+a_3\bar E, 
\]
where $a_3=2\beta_1+\beta_3-\alpha \ge 9/8$.
Since $\dim |3A|>3$, $\bar S_3$ is $\bar f$-horizontal.
Hence, $a_3=2$ and $\OOO_{\bar S_1}(\bar S_3)=\OOO_{\bar S_1}(\bar E)=\OOO_{\PP^2}(1)$. 
From the exact sequence 
\[
0 \longrightarrow 
\OOO_{\bar X}(\bar S_3-\bar S_1)
\longrightarrow
\OOO_{\bar X}(\bar S_3)
\longrightarrow
\OOO_{\bar S_1}(\bar S_3)
\longrightarrow 0.
\]
we get $\dim H^0(\OOO_{\bar X}(\bar S_3-\bar S_1))\ge 
9-3=6$. Therefore, $\dim |2A|=\dim |S_3-S_1|\ge 5$, a contradiction.
\end{case}

\begin{case} \textbf{Case \ref{no-q=5-B237}.}
We have
$\beta_1\ge 5\alpha$ and
$a_1\ge 14\alpha\ge 2$.
Moreover, $f(E)\in X$ is a cyclic quotient singularity of index $r=2$, $3$ or $7$.
If $r=2$, then $a_1\ge 7$. Hence the contraction $\bar f$ is birational
and $\hat q =5s_1+a_1e\ge 13$, $s_1\ge 2$, a contradiction.
If $r=3$, then $\beta_1=2/3+m_1$ and
$a_1=3+5m_1$, where $m_1\ge 1$ because $\beta_1\ge 3\alpha$. In this situation, $\bar f$ is birational and
\[
\hat q =5s_1+a_1e=5s_1+ 3e+5m_1e\ge 13. 
\]
On the other hand, $\g(\hat X)\ge 25$. This contradicts Theorem
\ref{theorem-main-10-19}.

Therefore, $r=7$.
Then 
$\beta_1=3/7+m_1$, $\beta_2=6/7+m_2$, $a_1=2+5m_1$, and 
$a_2=2+m_1+2m_2$.
If $\bar f$ is not birational, then $a_1=a_2=2$ and the divisors
$\bar S_1$, $\bar S_2$ are $\bar f$-vertical.
Since $\dim |\Theta|=\dim |\bar S_1|=1$, $\hat X$ is either a curve
or $\PP(1^2,2)$ (see \ref{case-not-birational-del-Pezzo}). 
On the other hand, $\dim |2\Theta|=\dim |\bar S_2|=4$,
a contradiction.
Hence the contraction $\bar f$ is birational. In this case,
\[
\hat q=2e+5(s_1+m_1e)=2e+ s_1+m_1e+2(s_2 +m_2e),
\]
where
$0<s_2 +m_2e=2(s_1+m_1e)$. So, $\hat q\ge 7$.
Since $\g (\hat X)\ge \g(X)=25$, by 
Theorem \ref{theorem-main-10-19} and of \ref{theorem-main-q=9}-\ref{theorem-main-q=7} of Theorem \ref{theorem-main}
we have $\hat X\simeq \PP(1^2,2,3)$.
Therefore, 
$s_2\le 2$ and $3=\dim |2\Theta|\ge \dim |S_2|=4$, a contradiction.
\end{case}

\begin{case} \textbf{Case \ref{no-q=5-B79}.}
As above we see that $f(E)$ is a non-Gorenstein point.
If $f(E)$ is the point of index $9$, then $\beta_1=2/9+m_1$ and
$a_1=1+5m_1$, where $m_1\ge 1$. Hence the contraction $\bar f$ is birational.
Then $\hat q =5s_1+ e+5m_1e\ge 11$, $s_1\ge 2$, 
$\hat q\ge 17$, $s_1\ge 5$, and $\hat q>19$ a contradiction.

Therefore, $f(E)$ is the point of index $7$.
Then 
$\beta_1=3/7+m_1$, $\beta_2=6/7+m_2$, $a_1=2+5m_1$, and
$a_2=2+m_1+2m_2$.
If $\bar f$ is not birational, then $a_1=a_2=2$ and divisors
$\bar S_1$, $\bar S_2$ are $\bar f$-vertical.
Since $\dim |\bar S_2|>2$, the variety $\hat X$ is not a curve.

\begin{scase}\label{no-q=5-B79-s-case-non-birational}
Therefore, $\hat X\simeq \PP(1,1,2)$ and the linear system $|\bar S_2|$
is base point free. In this situation, $\bar S_2$ is a Hirzebruch surface 
$\FF_n$. Let $\Lambda$ and $\Sigma$ be its fiber and negative section, respectively.
\[
-K_{\bar S_2}=(-K_{\bar X}-2\bar S_1)|_{\bar S_2} =
(3\bar S_1+2\bar E)|_{\bar S_2}
\]
Thus,
\[
2 \bar E|_{\bar S_2}=-K_{\bar S_2}-3\Lambda=
2\Sigma +(n-1)\Lambda.
\]
Since $\bar E|_{\bar S_2}$ is an irreducible section, 
there is only one possibility: $\bar E|_{\bar S_2}=\Sigma$ and $n=1$.
Since $\dim |\hat S_4|> 9$, the divisor $\bar S_4$ is $\bar f$-horizontal.
Similar to \eqref{equation-rel-Cl-79-1} and \eqref{equation-rel-Clb-1-79}
we have 
\begin{equation*}
-K_{\bar X} \sim \bar S_1+\bar S_4+a_4\bar E,
\qquad 
a_4=\beta_1+\beta_4-\alpha>0.
\end{equation*}
Hence, $a_4=1$ and $\bar S_4\sim \bar E+4\bar S_1$. 
From the exact sequence 
\[
0 \longrightarrow 
\OOO_{\bar X}(\bar S_4-\bar S_2)
\longrightarrow
\OOO_{\bar X}(\bar S_4)
\longrightarrow
\OOO_{\bar S_2}(\bar S_4)
\longrightarrow 0.
\]
we get $\dim H^0(\OOO_{\bar X}(\bar S_4-\bar S_2))\ge 
14-9=5$. Therefore, 
\[
\dim |2A|=\dim |S_4-S_2|\ge 4, 
\]
a contradiction.
\end{scase}
\begin{scase}\label{no-q=5-B79-s-case-birational}
Hence the contraction $\bar f$ is birational. Then
\[
\hat q=2e+5(s_1+m_1e)=2e+ s_1+m_1e+2(s_2 +m_2e),
\]
$0<s_2 +m_2e=2(s_1+m_1e)$, and $\hat q\ge 7$.
If $s_1+m_1e>1$, then $s_1+m_1e=3$, $\hat q=2e+15$,
$e\ge 2$, $\hat q\ge 19$, $e\ge 3$, and $s_1=0$, a contradiction.
Therefore, $s_1+m_1e=1$, $s_2 +m_2e=2$, $s_1=1$, $\hat q=2e+5$,
$\hat q=7$, and $e=1$. By \cite[Th. 1.4 (vi)]{Prokhorov2008a} 
$\hat X\simeq \PP(1^2,2,3)$. Further,
$s_2=2$ and $\dim |2\Theta|= \dim |S_2|$.
Using \eqref{equation-discrepancies-second} write
\[
\begin{array}{rll}
K_{\bar X}+\frac72\bar S_2&=&\bar f^*(K_{\hat X}+\frac72\hat S_2)
+(b-\frac72\gamma_2)\bar F.
\end{array}
\]
Pushing this relation down to $X$ we get
$2= b-7\gamma_2/2$.
Since $b\ge 2$, $Q:=\bar f(\bar F)\in \PP(1^2,2,3)$ is a smooth point \cite{Kawamata-1996}.
On the other hand, the base locus of $|2\Theta|$ is the point
$(0:0:0:1)$ of index $3$.
Hence, $\gamma_2=0$ and $b=2$.
Thus $\bar f$ is the usual blowup of a smooth point \cite{Kawakita-2001}.

Let $\hat {\SSS}_2$ be the subsystem of $|2\Theta|=|\hat S_2|$
consisting of divisors passing through $Q$.
Write
\[
\bar S_2\sim \bar f^* \hat {\SSS}_2=\bar {\SSS}_2 +\bar F,
\]
where $\bar {\SSS}_2$ is the birational transform.
Pushing these relations down to $X$ we get
$S_2\sim {\SSS}_2 + F$. 
Hence, $S_1\sim {\SSS}_2$. On the other hand,
$\dim {\SSS}_2=2$, a contradiction.
\end{scase}\end{case}

\begin{case} \textbf{Case \ref{no-q=5-B47}.}
Then $\beta_1\ge 3\alpha$, $a_1=5\beta_1-\alpha\ge 14\alpha\ge 2$.
Moreover, $f(E)$ is a non-Gorenstein point.
If $f(E)$ is the point of index $4$, then $\beta_1=1/4+m_1$,
$a_1=1+5m_1$, where $m_1\ge 1$. Hence the contraction $\bar f$ is birational. Then
$\hat q =5s_1+a_1e=5s_1+ e+5m_1e\ge 11$, $s_1\ge 2$, 
$\hat q\ge 17$, $s_1\ge 5$, and 
$\hat q>19$, a contradiction.

Therefore, $f(E)$ is the point of index $7$.
Then 
$\beta_1=3/7+m_1$, $\beta_2=6/7+m_2$, $a_1=2+5m_1$, and
$a_2=2+m_1+2m_2$.
If $\bar f$ is not birational, then $a_1=a_2=2$,
$\bar S_1$, $\bar S_2$ are $\bar f$-vertical.
Since $\dim |\bar S_2|>1$, $\hat X$ is not a curve.
Hence, $\hat X\simeq \PP(1,1,2)$. 
Then we get a contradiction as in 
\ref{no-q=5-B79-s-case-non-birational}.

Hence the contraction $\bar f$ is birational. Then
\[
\hat q=2e+5(s_1+m_1e)=2e+ s_1+m_1e+2(s_2 +m_2e).
\]
This gives us
$0<s_2 +m_2e=2(s_1+m_1e)$ and $\hat q\ge 7$.
If $s_1+m_1e>1$, then $s_1+m_1e=3$, $\hat q=2e+15$,
$e\ge 2$, $\hat q\ge 19$, $e\ge 3$, and $\hat q> 19$, a contradiction.
Therefore, $s_1+m_1e=1$, $s_2 +m_2e=2$, $s_1=1$, $\hat q=2e+5$,
$\hat q=7$, $e=1$, and $\hat X\simeq \PP(1^2,2,3)$. 
The we get a contradiction as in 
\ref{no-q=5-B79-s-case-birational}.
\end{case}

\begin{case} \textbf{Case \ref{no-q=5-B412}.}
Then $\beta_1\ge 5\alpha$, $a_1=5\beta_1-\alpha\ge 24\alpha\ge 2$.
In particular, this implies that $f(E)$ is a non-Gorenstein point.
Assume that $f(E)$ is the point of index $4$.
Near $f(E)$ we have $-K_X\sim A$, so $\beta_1=1/4+m_1$, where $m_1\ge 1$.
Hence, $a_1=5\beta_1-\alpha=1+5m_1\ge 6$.
This implies that $\bar f$ is birational, $\hat q=5s_1+e+5m_1e\ge 11$,
$\Cl(X)\simeq \ZZ$, $s_1\ge 2$, $\hat q\ge 17$, $s_1\ge 5$, and $\hat q>19$,
a contradiction.
Therefore, $f(E)$ is the point of index $12$.
Then we can argue as in the case \ref{no-q=5-B47}
\end{case}

\begin{case} \textbf{Case \ref{no-q=5-B223}.}
We show that in this case $X$ is a hypersurface of degree $4$ in $\PP(1^2,2^2,3)$.
The following is similar to Proposition \ref{lemma-DuVal-5}.
\begin{stheorem}{\bf Lemma.}
Let $M$ be a del Pezzo surface with Du Val singularities such that $K_M^2=6$.
The following are equivalent:
\begin{enumerate}
 \item 
$-K_M\qq 3L$ for some Weil divisor $L$;
 \item 
$\Sing(M)$ consists of 
a unique Du Val point of type $A_2$.
 \item 
$M\simeq M_4\subset \PP(1^2,2,3)$, where $M_4$
is given by $x_1x_3+x_2^2+x_1'^4=0$. 
\end{enumerate}
\end{stheorem}

Similar to \ref{explanation-DP} it is sufficient to prove the following. 

\begin{stheorem}{\bf Proposition.}
Let $X$ be a Fano threefold of type \textup{\ref{no-q=5-B223}}.
Then a general member $M\in |2A|$ is a Du Val del Pezzo surface
of degree $6$ such that $-K_M$ is divisible by $3$ in $\Cl(M)$. 
\end{stheorem}

\begin{proof}
Put $\MMM:=|2A|$.
As in \ref{proposition-q=9-D5} we show that the pair $(X,\MMM)$ is 
canonical.
Assume the converse.
In \eqref{equation-discrepancies} we have 
$0<a_0=\beta_2-\alpha$. Hence, $\beta_2>\alpha$ and $\beta_1> \frac12 \alpha$.
By \eqref{equation-rel-Clb-1-79} we have 
\[
a_1=5\beta_1-\alpha>3\beta_1,\quad a_2=\beta_1+2\beta_2-\alpha>\beta_1+\beta_2,
\qquad a_1,\, a_2\ge 2.
\]
If $\bar f$ is not birational, then $\bar S_1$ and $\bar S_2$ are $\bar f$-vertical.
Since $\dim |\bar S_2|=4$, $\hat X$ is not a curve.
Since $\dim |\bar S_1|=1$, the surface $\hat X$ is either $\PP^2$ or $\PP(1^2,2)$.
We have a contradiction because $\dim |\bar S_2|=4$.
Therefore, the contraction $\bar f$ is birational. Write
\[
\hat q=5s_1+a_1e=s_1+2s_2+a_2e\ge 7.
\]
If $s_1\ge 2$, then $\hat q\ge 13$ and $\alpha \ge 1$ because $\g(\hat X)\ge 21$.
In this situation, 
we get successively $\beta_1>1$, $a_1\ge 4$, $\hat q\ge 17$, $s_1\ge 5$, and $\hat q>19$,
a contradiction. Therefore, $s_1=1$. By \cite[Th. 1.4 (vi)]{Prokhorov2008a}
$\hat X\simeq \PP(1^2,2,3)$. Then $3=\dim |2\Theta|>\dim |2A|=4$.
Hence, $s_2\ge 3$ and $a_2=0$, a contradiction.
\end{proof}
\end{case}

\begin{case} \textbf{Cases \ref{no-q=5-B222336}.}
Take $\MMM=|2A|$. Near the point of index $6$ we have 
$\MMM\sim -4K_X$. Hence, $c\le 1/4$, $\beta_2\ge 4\alpha$, and $\beta_1\ge 2\alpha$.
If $\alpha\ge 1$, then $a_2\ge 9\alpha\ge 9$ by \eqref{equation-rel-Clb-1-79}.
Hence the contraction $\bar f$ is birational and $\hat q\ge 2s_2+9\ge 11$.
In this case, $s_2\ge 4$, $\hat q\ge 17$, $s_2>7$, and $\hat q>19$,
a contradiction. Therefore, $\alpha<1$ and $f(E)$ is a point
of index $r=2$, $3$, or $6$.

\begin{scase} \textbf{Subcase $r=2$.}
 Then $\alpha=1/2$, $\beta_1=1/2+m_1$, and $\beta_2=m_2$, where $m_i\in \ZZ$,
$m_1\ge 1$, $m_2\ge 2$.
Hence, $a_1=2+5m_1\ge 7$ and $a_2=2m_2+m_1$.
This implies that $\bar f$ is birational and
\[
\hat q= 2e+5(s_1+em_1)=2(s_2+em_2)+s_1+em_1\ge 7.
\]
If $s_1+em_1>1$, then $\hat q\ge 13$ and so $s_1+em_1$ is odd.
Hence, $\hat q\ge 17$, $e\ge 2$, $\hat q\ge 19$, $e\ge 3$, and $\hat q>19$,
a contradiction. Thus, $s_1+em_1=1$ and so 
$m_1=e=1$, $s_1=0$, $\hat q= 7$, $s_2+em_2=3$, $s_2=1$,
This contradicts \cite[Th. 1.4 (vi)]{Prokhorov2008a}.
\end{scase}

\begin{scase} \textbf{Subcase $r=3$ and $\alpha=2/3$.}
Then $\beta_1=1/3+m_1$ and $\beta_2=2/3+m_2$, where $m_i\in \ZZ$,
$m_1\ge 1$, and $m_2\ge 2$.
Hence, $a_1=1+5m_1\ge 6$ and $a_2=1+2m_2+m_1$.
This implies that $\bar f$ is birational and
\[
\hat q= e+5(s_1+em_1)=e+2(s_2+em_2)+s_1+em_1\ge 7.
\]
By \cite[Th. 1.4 (vi)]{Prokhorov2008a} $s_2\ge 2$, so 
$\hat q\ge 9$, $\Cl(\hat X)$ is torsion free, 
 $s_2\ge 4$, $\hat q\ge 17$, $s_2>7$,
and $\hat q>19$, a contradiction. 
\end{scase}

\begin{scase} \textbf{Subcase $r=3$ and $\alpha=1/3$.}
Then $\beta_1=2/3+m_1$ and $\beta_2=1/3+m_2$, where $m_i\in \ZZ$,
$m_2\ge 1$.
Hence, $a_1=3+5m_1$ and $a_2=1+2m_2+m_1\ge 3$.
If $\bar f$ is not birational, then it is generically 
$\PP^2$-bundle and the divisors $\bar S_1$ and $\bar S_2$
are $\bar f$-vertical. Since $\dim |\bar S_2|=3$,
$\bar S_2$ must be reducible, a contradiction.
Therefore, the contraction $\bar f$ is birational and
\[
\begin{array}{l}
\hat q=3e+5(s_1+em_1)=e+2(s_2+em_2)+s_1+em_1\ge 5,
\\[5pt]
e+2(s_1+em_1)=s_2+em_2>em_2.
\end{array}
\]
It is easy to see from the second relation that $s_1+em_1\neq 0$.
If $s_1+em_1=1$, then $2=s_2+e(m_2-1)$, $s_2\le 2$, $\hat q\le 8$,
$e=1$, and $\hat q=8$. 
By Lemma \ref{lemma-dk} and Theorem \ref{theorem-main}
\ref{theorem-main-q=8} we get a contradiction.
Thus $s_1+em_1>1$. As above we have $\hat q=11$.
But then $e<0$, a contradiction.
\end{scase}

\begin{scase} \textbf{Subcase $r=6$.}
Then $\alpha=1/6$, $\beta_1=5/6+m_1$, and $\beta_2=4/6+m_2$, where $m_i\in \ZZ_{\ge 0}$.
Hence, $a_1=4+5m_1$ and $a_2=2+2m_2+m_1$.
This implies that $\bar f$ is birational and
\[
\hat q= 4e+5(s_1+em_1)=2e+2(s_2+em_2)+s_1+em_1.
\]
If $s_1+em_1>0$, then $\hat q\ge 9$, $\Cl(\hat X)$ is torsion free, 
$s_2\ge 4$, $\hat q\ge 11$. On the other hand, $s_2\le (\hat q-3)/3\le 8$.
By Theorem \ref{theorem-main-10-19} we get successively $\hat q \le 13$,
$s_2\le 5$, $s_1+em_1=1$, $\hat q= 4e+5$,
$s_2+em_2=e+2$, $e=2$, $m_2=0$, $s_2=4$, $\hat q=14$, a contradiction.
Therefore, $s_1=m_1=0$ and
$\hat q= 4e=2e+2(s_2+em_2)$.
By Lemma \ref{lemma-torsion-ae} the group $\Cl(\hat X)$ is torsion free 
and $e=1$.
Hence, $\hat q=4$ and $s_2+em_2=1$.
By Proposition \ref{proposition-dim-A}
$\hat X\simeq \PP^3$.
Using \eqref{equation-discrepancies-second} we write
\[
K_{\bar X}+4\bar \MMM=\bar f^*(K_{\hat X}+4\hat \MMM)+(b-4\gamma_2) \bar S_1,
\]
where $K_{\hat X}+4\hat \MMM\sim 0$. Hence, $b=4\gamma_2+3$
and so $Q:=\bar f(\bar S_1)$ is a point.
Since $\hat \MMM=|\OOO_{\PP^3}(1)|$ is a complete linear system and it is 
a base point free, we have $\gamma_2=0$ and $b=3$.
Let $\hat \MMM'\subset |\OOO_{\PP^3}(1)|$ be the 
subsystem consisting of divisors passing through $Q$.
Then 
\[
K_{\bar X}+4\bar \MMM'=\bar f^*(K_{\hat X}+4\hat \MMM')+(b-4\gamma_2')\bar S_1. 
\]
where $\gamma_2'>0$. 
Then $-8+4m=b- 4\gamma_2'$, where $\MMM'\sim mA$.
We obtain $8=3+4m+4\gamma'_2$, where $m>1$ because 
$\dim \MMM'>\dim |A|$, a contradiction.
\end{scase}
\end{case}

\begin{mtheorem}
\label{q=5-general-case}
{\bf Lemma.}
Let $X$ be a $\QQ$-Fano threefold with $\qQ(X)=5$.
If $\qW(X)\neq \qQ(X)$, then $\g(X)\le 8$. 
\end{mtheorem}

\begin{proof}
Similar to Lemma \ref{lemma-q=7-torsion-general}.
\end{proof}

\section{Case $\qQ(X)=4$}
\begin{mtheorem}{\bf Lemma.}
\label{lemma-q=4-table}
Let $X$ be a $\QQ$-Fano threefold with $\qW(X)=\QQ(X)=4$ and 
$\g(X)\ge 22$. Assume that $\B(X)\neq \emptyset$. Let
$-K_X=4A$.
Then the group $\Cl(X)$ is torsion free and we have one of the following cases:
{\rm
\setlongtables\renewcommand{\arraystretch}{1.3}
\begin{longtable}{|c|c|c|c|c|c|c|c|}
\hline
&&&\multicolumn{5}{c|}{$\dim |kA|$}
\\
\hhline{|~|~|~|-----}
&$\B$&$A^3$& $|A|$&$|2A|$&$|3A|$&$|-K|$
\\[5pt]
\hline
\endfirsthead
\hline
&$\B$&$A^3$& $|A|$&$|2A|$&$|3A|$&$|-K|$
\\[5pt]
\hline
\endhead
\hline
\endlastfoot
\hline
\endfoot

\nr\label{no-q=4-B11a}&
$(11)$&
$10/11$& $2$& $7$& $16$& $30$

\\
\nr\label{no-q=4-B57}&
$(5, 7)$&
$32/35$& $2$& $7$& $16$& $30$

\\
\nr\label{no-q=4-B35}&
$(3, 5)$&
$13/15$& $2$& $7$& $16$& $29$

\\
\nr\label{no-q=4-B39}&
$(3, 9)$&
$8/9$& $2$& $7$& $16$& $29$

\\ \nr\label{no-q=4-B5}&
$(5)$&
$4/5$& $2$& $7$& $15$& $27$

\\
\nr\label{no-q=4-B13}&
$(13)$&
$10/13$& $1$& $6$& $14$& $25$

\\
\nr\label{no-q=4-B59}&
$(5, 9)$&
$34/45$& $1$& $6$& $13$& $25$

\\
\nr\label{no-q=4-B7}&
$(7)$&
$5/7$& $2$& $6$& $13$& $24$

\\
\nr\label{no-q=4-B11b}&
$(11)$&
$8/11$& $2$& $6$& $13$& $24$

\\
\hline
&&&&&&

\\[-12pt]
\nr\label{no-q=4-B3}&
$(3)$& $2/3$&$2$&$6$&$13$&$23$
\end{longtable}
}
\end{mtheorem}

\begin{case}
Now let $X$ be a $\QQ$-Fano threefold with $\qQ(X)=4$ and $\g(X)\ge 22$.
We also assume that $\qQ(X)=\qW(X)$. For the case $\qQ(X)\neq \qW(X)$ see 
Lemma \ref{q=4-general-case} below. 
If $\B(X)=\emptyset$, then $X\simeq\PP^3$ (see, e.g., \cite{Iskovskikh-Prokhorov-1999}).
Hence,
$X$ is described in Lemma \ref{lemma-q=4-table} and $\Cl(X)$ is torsion free.
We will show that cases \ref{no-q=4-B11a} -- \ref{no-q=4-B11b} do not occur.
The case \ref{no-q=4-B3} is considered in \cite{Sano-1996}.

Thus we may assume that $\B(X)\neq \emptyset$ and $\g(X)\ge 23$.
Put $\MMM:=|A|$.
Let $P\in X$ be a point of index $r'$, where we take $r'$ as follows:
\[
\begin{array}{c|ccccccccc}
\text{No.}
&\text{\xref{no-q=4-B11a}}
&\text{\xref{no-q=4-B57}}
&\text{\xref{no-q=4-B35}}
&\text{\xref{no-q=4-B39}}
&\text{\xref{no-q=4-B5}}
&\text{\xref{no-q=4-B13}}
&\text{\xref{no-q=4-B59}}
&\text{\xref{no-q=4-B7}}
&\text{\xref{no-q=4-B11b}}
\\
r' & 11& 5 & 5 & 9 & 5 & 13 & 9 & 7 & 11
\\ 
m & 3& 4 & 4 & 7 & 4 & 10 & 7 & 2 & 3
\end{array}
\]
Here $m$ is an integer such that $0\le m<r'$
and $\MMM\sim -mK_X$ near $P$. 
Then $\alpha/\beta_1=c\le 1/m$ (see \ref{lemma-c-42}) and so 
$\beta_1\ge m\alpha$.
We have
\begin{equation}
\label{equation-KS-q=4}
\begin{array}{llll}
K_{\bar X}&+4\bar S_1&+a_1\bar E&\sim 0,
\\[5pt]
K_{\bar X}&+2\bar S_2&+a_2\bar E&\sim 0,
\end{array}
\end{equation}
where
\begin{equation}
\label{equation-q=4-a}
a_1=4\beta_1-\alpha\ge (4m-1)\alpha>0,
\quad
a_2=2\beta_2-\alpha.
\end{equation}
\end{case}

\begin{stheorem}
{\bf Claim.} 
 $f(E)$ is a point of index $>1$.
\end{stheorem}
\begin{proof}
Assume that $f(E)$ is either a curve or Gorenstein point.
Then $a_1\ge 7$, so
 the contraction $\bar f$ is birational and 
we get successively $\hat q\ge 4s_1+7\ge 11$,
$s_1\ge 2$, $\hat q\ge 17$, $s_1\ge 5$, and $\hat q>19$, a contradiction.
Therefore, $f(E)$ is a cyclic quotient singularity of index $r>1$ and $\alpha=1/r$.
We also claim that $a_2>0$. Indeed, otherwise $2\beta_2=\alpha=1/r$.
On the other hand, $\beta_2\in \frac 1r \ZZ$, a contradiction.
\end{proof}

\begin{case}
Assume that $\bar f$ is not birational. 
By \eqref{equation-KS-q=4} the divisor $\bar S_1$ is $\bar f$-vertical and 
the class of the divisor $\bar f(\bar S_1)$ is a generator of $\Cl(\hat X)$.
If $\hat X$ is a curve, then $\dim |\bar S_1|=1$
(cases \xref{no-q=4-B13} and \xref{no-q=4-B59}).
Since $\dim |\bar S_2|>2$, $\bar S_2$ is $\bar f$-horizontal.
Then $\bar f$ is a generically $\PP^2$-bundle,
$a_2=1$, and $\OOO_{\bar S_1}(\bar S_2)=\OOO_{\bar S_1}(\bar E)=\OOO_{\PP^2}(1)$. 
Thus $\bar S_1\simeq \PP^2$ is a general fiber.
From the exact sequence 
\[
0 \longrightarrow 
\OOO_{\bar X}(\bar S_2-\bar S_1)
\longrightarrow
\OOO_{\bar X}(\bar S_2)
\longrightarrow
\OOO_{\bar S_1}(\bar S_2)
\longrightarrow 0.
\]
we get $\dim H^0(\OOO_{\bar X}(\bar S_2-\bar S_1))\ge 
7-3=4$. Therefore, $\dim |A|=\dim |S_2-S_1|\ge 3$, a contradiction.

Thus $\hat X$ is a surface. In cases \xref{no-q=4-B13} and \xref{no-q=4-B59}
we have $a_1\ge 3$ by \eqref{equation-q=4-a}. Hence these cases are impossible
and so $\dim |\bar S_1|=2$.
Then $\hat X\simeq \PP^2$. Since $\dim |\bar S_2|>5$, $\hat S_2$ is $\bar f$-horizontal.
Thus $a_2=0$, a contradiction.
\end{case}
\begin{case}
Therefore, the contraction $\bar f$ is birational. 
Then
\[
\hat q =4s_1+a_1e\ge 5, \quad \g(\hat X)\ge \g (X)\ge 23. 
\]
Hence, $\Cl(\hat X)$ is torsion free. 
Consider the case $\hat q\ge 6$. Then
$\hat X\simeq \PP(1^2,2,3)$ by \ref{theorem-main-q=9}-\ref{theorem-main-q=6} of
Theorem \ref{theorem-main}. Moreover, $s_1=1$ and 
$\dim |S_1|\le \dim |\Theta|=1$.
Thus we are in cases \xref{no-q=4-B13} or \xref{no-q=4-B59}.
By \eqref{equation-q=4-a} we have $a_1\ge 3$, so $e=1$.
Further, $\dim |\hat S_2|\ge 6$, so $s_2=3$, $a_2=1$, and $|S_2|$
is the birational transform of $|3\Theta|$. 
If $\bar f(\bar F)$ is a non-Gorenstein point, then
$\bar f$ is a Kawamata blowup \cite{Kawamata-1996}
and our link is toric. This contradicts \cite[\S 10]{Prokhorov2008a}.
Thus $\bar f(\bar F)$ is either a curve or a smooth point.
On the other hand, the only base point of the 
linear system $|\hat S_2|=|3\Theta|$
is the point of type $\frac12(1,1,1)$.
Hence, $\bar S_2=\bar f^*\hat S_2$.
Since $e=1$ and $\Cl(\hat X)$ is torsion free, by Lemma \ref{lemma-torsion-ae}
we have $F\sim A$.
Using \eqref{equation-discrepancies-second} write 
\[
K_{\bar X}+\frac73\bar S_2=\bar f^*\Bigl(K_{\hat X}+\frac73\hat S_2\Bigr)+b\bar F.
\]
Pushing this relation down to $X$ we get $b=-4+2\cdot 7/3=2/3\notin \ZZ$, a contradiction.

Therefore, $\hat q=5$ and $s_1=a_1=e=1$. Then 
by \ref{theorem-main-q=5} of
Theorem \ref{theorem-main} we have either $\hat X\simeq \PP(1^3,2)$ 
or $\g(\hat X)=23$ (cases \xref{no-q=4-B7} or \ref{no-q=4-B11b}).
In the latter case, $\dim |\Theta|\ge \dim |S_1|= 2$, so
again $\hat X\simeq \PP(1^3,2)$ by Proposition \ref{proposition-dim-A}.
Since $a_1=1$, by \eqref{equation-q=4-a} we infer
$m \le 3$. Hence, we are in cases \ref{no-q=4-B11a}, \ref{no-q=4-B7} or \ref{no-q=4-B11b}.
As above we have $F\sim A$ and 
\[
K_{\bar X}+5\bar S_1=\bar f^*(K_{\hat X}+5\hat S_1)+(b-5\gamma_1)\bar F.
\]
Thus $b=1+5\gamma_1\ge 1$.
Since $\dim |\hat S_1|= \dim |S_1|$, 
we have $\gamma_1=0$ and $b=1$. Since $\PP(1^3,2)$ has no 
Gorenstein singular points, $C:=\bar f(\bar F)$ is a curve.
Clearly, $C$ is contained in $\hat E$, where $\OOO_{\hat X}(\hat E)=\OOO_{\PP(1^3,2)}(1)$.
By \cite{Kawamata-1996} the curve $C$ does not pass through the singular point of $\PP(1^3,2)$.
This implies that $C$ is a Cartier divisor on $\hat E\simeq \PP(1^2,2)$
and so $C$ is a complete intersection of members $\hat E\in |\OOO_{\PP(1^3,2)}(1)|$
and $\hat D\in |\OOO_{\PP(1^3,2)}(a)|$, $a>0$.
By \eqref{equation-KS-q=4} $s_2\le 2$.
On the other hand, $\dim |\hat S_2|\ge \dim |S_2|= 6$.
Hence, the case \ref{no-q=4-B11a} does not occur, $s_2=2$, and 
$|2\Theta|$ is the birational transform of $|2A|$.
As above, using \eqref{equation-discrepancies-second} we obtain
$\gamma_2=0$ and $\delta =1$. Hence $\hat E$ is the only element 
of $|\OOO_{\PP(1^3,2)}(1)|$ passing through $C$
and no reduced irreducible members $|\OOO_{\PP(1^3,2)}(2)|$ pass through $C$.
Therefore, $a\ge 3$. In this case, $\bar E$ is covered by a family of curves 
having negative intersection number with $-K_{\bar X}$, a contradiction.
\end{case}

\begin{mtheorem}
\label{q=4-general-case}
{\bf Lemma.}
Let $X$ be a $\QQ$-Fano threefold with $\qQ(X)=4$.
Assume that the torsion part of $\Cl(X)$ is non-trivial.
Then $\g(X)\le 19$. 
\end{mtheorem}

\section{Case $\qQ(X)=3$}
\begin{mtheorem}{\bf Lemma.}\label{lemma-q=3-table}
Let $X$ be a $\QQ$-Fano threefold with $\qW(X)=\qQ(X)=3$ and 
$\g(X)\ge 21$. Assume that $\B(X)\neq \emptyset$. Let
$-K_X=3A$.
Then the group $\Cl(X)$ is torsion free and we have one of the following cases:

\rm
\par\medskip\noindent
\setlongtables\renewcommand{\arraystretch}{1.3}
\begin{longtable}{|c|c|c|c|c|c|c|c|}
\hline
&&&\multicolumn{5}{c|}{$\dim |kA|$}
\\
\hhline{|~|~|~|-----}
&$\B$&$A^3$& $|A|$&$|2A|$&$|-K|$
\\[5pt]
\hline
\endfirsthead
\hline
&$\B$&$A^3$& $|A|$&$|2A|$&$|-K|$
\\[5pt]
\hline
\endhead
\hline
\endlastfoot
\hline
\endfoot

\nr\label{no-q=3-B4}& 
$(4)$&
$9/4$& $4$& $14$& $32$

\\ \nr\label{no-q=3-B7}& 
$(7)$&
$16/7$& $4$& $14$& $32$

\\ \nr\label{no-q=3-B244}& 
$(2, 4, 4)$&
$2$& $3$& $12$& $28$

\\ \nr\label{no-q=3-B45}& 
$(4, 5)$&
$37/20$& $3$& $11$& $26$

\\ \nr\label{no-q=3-B24}& 
$(2, 4)$&
$7/4$& $3$& $11$& $25$

\\ \nr\label{no-q=3-B27}& 
$(2, 7)$&
$25/14$& $3$& $11$& $25$

\\ \nr\label{no-q=3-B210}& 
$(2, 10)$&
$9/5$& $3$& $11$& $25$

\\ \nr\label{no-q=3-B5}& 
$(5)$&
$8/5$& $3$& $10$& $23$

\\ \nr\label{no-q=3-B8}& 
$(8)$&
$13/8$& $3$& $10$&$23$

\\ \nr\label{no-q=3-B11}& 
$(11)$&
$18/11$& $3$& $10$&$23$

\\ \nr\label{no-q=3-B227}& 
$(2, 2, 2, 7)$&
$23/14$& $2$& $10$& $23$

\\ \nr\label{no-q=3-B57}& 
$(5, 7)$&
$54/35$& $2$& $9$& $22$

\\
\hline
&&&&&
\\[-9pt] \nr\label{no-q=3-B2}& 
$(2)$&
$3/2$& $3$& $10$& $22$

\end{longtable}
\end{mtheorem}

\begin{case}
Now let $X$ be a $\QQ$-Fano threefold with $\qQ(X)=\qW(X)=3$ and $\g(X)\ge 21$.
If $\B(X)=\emptyset$, then $X$ is a smooth quadric in $\PP^4$ (see, e.g., \cite{Iskovskikh-Prokhorov-1999}).
Thus we may assume that $\B(X)\neq \emptyset$.
Then
$X$ is described in Lemma \ref{lemma-q=3-table} and $\Cl(X)$ is torsion free.
In the case \ref{no-q=3-B2} we have 
$X\simeq X_3\subset \PP(1^4,2)$ by \cite{Sano-1996}.
We will show that cases
\ref{no-q=3-B4} -- \ref{no-q=3-B57} do not occur.

Write 
\begin{equation}
\label{equation-q=3-K}
K_{\bar X}+3\bar S_1 +a_1\bar E\sim 0, 
\quad 
K_{\bar X}+\bar S_1+\bar S_2 +a_2\bar E\sim 0, 
\end{equation}
where
\begin{equation}
\label{equation-q=3-aa}
a_1=3\beta_1-\alpha,
\quad
a_2=\beta_1+\beta_2-\alpha. 
\end{equation}
In all cases we have $\dim |A|\ge 2$. Put $\MMM:=|A|$.
\end{case}
\begin{case} {\bf Cases \ref{no-q=3-B4}-\ref{no-q=3-B11}.}
Let $P\in X$ be a point of index $r'$, where we take $r'$ as follows:
\[
\begin{array}{c|cccccccccc}
\text{No.}

&\text{\xref{no-q=3-B4}}
&\text{\xref{no-q=3-B7}}
&\text{\xref{no-q=3-B244}}
&\text{\xref{no-q=3-B45}} 
&\text{\xref{no-q=3-B24}}
&\text{\xref{no-q=3-B27}}
&\text{\xref{no-q=3-B210}}
&\text{\xref{no-q=3-B5}}
&\text{\xref{no-q=3-B8}}
&\text{\xref{no-q=3-B11}}
\\
r' & 4 & 7 & 4 & 4 & 4 & 7 & 10 & 5 & 8 & 11
\\ 
m & 3 & 5 & 3 & 3 & 3 & 5 & 7 & 2 & 3 & 4
\end{array}
\]
Here $m$ is an integer such that $0\le m<r'$
and $\MMM\sim -mK_X$ near $P$. 
Then $\alpha/\beta_1=c\le 1/m$ (see \ref{lemma-c-42}) and 
so $\beta_1\ge m\alpha$. Hence,
\begin{equation}
\label{equation-q=3-a}
a_1=3\beta_1-\alpha\ge (3m-1)\alpha\ge 5\alpha>0. 
\end{equation}
\begin{stheorem}
{\bf Lemma.}
$\hat X\simeq \PP(1^2,2,3)$. 
\end{stheorem}

\begin{proof}
If $\bar f$ is not birational, then 
$\bar S_1$ is $\bar f$-horizontal (because $\dim |\bar S_1|\ge 3$).
Restricting \eqref{equation-q=3-K} to a general fiber 
we obtain $a_1=0$, a contradiction.

Therefore, the contraction $\bar f$ is birational. In this case, $\hat q\ge 3s_1+a_1\ge 4$.
If $f(E)$ is either a curve or Gorenstein point, then we get successively
$a_1\ge 5$, $\hat q\ge 8$, $s_1\ge 2$
(see \cite[Theorem 1.4 (vi)]{Prokhorov2008a}), $\hat q\ge 11$, $s_1\ge 4$, $\hat q\ge 17$,
$s_1\ge 7$, 
$\hat q> 19$, a contradiction.

Thus $f(E)$ is a point of index $>1$.
By Lemma \ref{lemma-discrepancies} $\alpha<1$.
Then by Lemma
\ref{lemma-dk}
$\g(\hat X)\ge \g(X)\ge 22$.
Using Theorem \ref{theorem-main-10-19} 
and \ref{theorem-main-q=9}-\ref{theorem-main-q=6} of Theorem \ref{theorem-main}
we get either
$\hat X\simeq \PP(1,2,3,5)$, $\PP(1^2,2,3)$, or $\hat q\le 5$.
If $\hat X\simeq \PP(1,2,3,5)$, then 
$\dim |3\Theta|=2$, $s_1\ge 4$, $\hat q\ge 3s_1\ge 13$, a contradiction.
Thus we may assume that $4\le \hat q\le 5$. 
Then by Lemmas 
\ref{q=5-general-case} and \ref{q=4-general-case}
the group $\Cl(\hat X)$ 
is torsion free. 
If $\hat X$ is singular, then by Proposition \ref{proposition-dim-A} we have 
$\dim |\Theta|\le 2$, so $s_1\ge 2$ and $\hat q\ge 7$, a contradiction.
It remains to consider the case where $\hat X$ is smooth.
Then $\hat q=4$ and $\hat X\simeq \PP^3$. Since $|2\Theta|=9$ and $\dim |2A|\ge 10$, 
we have $s_2\ge 3$.
On the other hand, 
$a_2=\beta_1+\beta_2-\alpha\ge (m-1)\alpha$.
Hence, $a_2\ge 1$. We get $4= \hat q\ge s_1+s_2+a_2\ge 5$.
The contradiction proves the lemma.
\end{proof}
\begin{scase}\label{case-q=3-933}
By the above, $s_1=2$ and $a_1=e=1$.
Hence, $\dim |2\Theta|=\dim |A|=3$ and the linear system $|2\Theta|$ is the 
proper transform of $|A|$. By Lemma \ref{lemma-torsion-ae} $F\sim A$.
Using \eqref{equation-discrepancies-second} write 
\[
\begin{array}{ll}
K_{\bar X}+\frac72\bar S_1\equiv
\bar f^*\left(K_{\hat X}+\frac72\hat S_1\right)+\left(b-\frac72\gamma_1\right)\bar F
\equiv\left(b-\frac72\gamma_1\right)\bar F,
\\[5pt]
K_{\bar X}+7\bar E\equiv\bar f^*\left(K_{\hat X}+7\hat E\right)+\left(b-7\delta\right)\bar F
\equiv\left(b-7\delta\right)\bar F.
\end{array}
\]
Pushing these relations down to $X$ we get
\[
 b=\frac12+\frac72\gamma _1=-3+7\delta.
\]
Hence, $b\ge 1/2$.
On the other hand, $\Bs|2\Theta|$ is the point of index $3$ and it is
a cyclic quotient singularity. By \cite{Kawamata-1996}
$\bar f(\bar F)\neq \Bs|2\Theta|$, so $\gamma _1=0$ and $b=1/2=\delta$.
Write
\[
K_{\hat X}+\hat E +6\hat B\sim 0, 
\]
where $\hat B$ is a general member 
of $|\Theta|$.
By taking the pull-back we get
\[
0\sim \bar f^*(K_{\hat X}+\hat E +6\hat B)\sim 
K_{\bar X}+\bar E +6\bar B+ \lambda \bar F,
\]
where $\lambda\ge \delta -b= 0$ and $\dim |\bar B|\ge 1$. Therefore,
\[
 K_{X}+6B+ \lambda F\sim 0.
\]
This contradicts $\qQ(X)=3$. 
\end{scase}\end{case}

\begin{case} \textbf{Cases  \ref{no-q=3-B227} and \ref{no-q=3-B57}.}
Near the point of index $7$ we have $-5K_X\sim A$.
By Lemma \ref{lemma-c-42} $c\le 1/5$.
Hence, $\beta_1\ge 5\alpha$ and $a_1\ge 14\alpha\ge 2$.

\begin{stheorem}
{\bf Lemma.}
$f(E)$ is a point of index $7$. 
\end{stheorem}
\begin{proof}
Assume the converse. 
If $f(E)$ is either a 
curve or Gorenstein point, then 
$\hat q\ge 3s_1+14e\ge 17$. Hence, $e=1$.
On the other hand, $|\Theta|=\emptyset$, a contradiction.
Therefore
$f(E)$ is a point of index $r=2$ (resp. $r=5$)
in the case \ref{no-q=3-B227} (resp. in the case \ref{no-q=3-B57}).
We have $a_1\ge 3$.
If $\hat f$ is not birational, then $\hat X\simeq\PP^1$ and $\bar f$ is a generically $\PP^2$-bundle.
Moreover, $\hat S_1$ is $\bar f$-vertical. This contradicts $\dim |S_1|=2$.
Therefore, $\hat f$ is birational. Then $\hat q= 3s_1+a_1e\ge 6$.
On the other hand, $\g (\hat X)\ge 21$, so 
$\hat X\simeq \PP(1,2,3,5)$, $\PP(1^2,2,3)$ or $ X_6\subset \PP(1^2,2,3,5)$ by Theorem \ref{theorem-main-10-19}
and \ref{theorem-main-q=9}-\ref{theorem-main-q=6} of Theorem \ref{theorem-main}.
Since $\dim |\hat S_1|\ge \dim |A|=2$, we have $|\hat S_1|\neq |\Theta|$, so $s_1\ge 2$,
$\hat q\ge 9$, 
and $\hat X\simeq \PP(1,2,3,5)$.
In this case, $s_1\ge 3$
and $\hat q \ge 16$, a contradiction.
\end{proof}

\begin{scase}
Thus $\alpha=1/7$, 
$\beta_1=5/7+m_1$, and $\beta_2=3/7+m_2$, where $m_1$, $m_2$ are 
non-negative integers. Hence, $a_1=2+3m_1\ge 2$ and $a_2=1+m_1+m_2\ge 1$.
\end{scase}

\begin{scase}
If $\bar f$ is not birational, then $m_1=0$, $a_1=2$, $\bar S_1$ is $\bar f$-vertical.
Since $\dim |\bar S_1|=2$, $\hat X$ is not a curve.
Thus $\bar f$ is a $\QQ$-conic bundle, $\hat X\simeq \PP^2$, $\bar f_*|\bar S_1|$
is the complete linear system of lines on $\hat X\simeq \PP^2$, and
$E$ is a generically section of $\bar f$. 
In this situation, there is an open subset $U\subset \hat X$
such that $\hat X\setminus U$ is finite, 
$\bar f^{-1}(U)$ is smooth, 
$\bar f$ induces a $\PP^1$-bundle $\bar f^{-1}(U)\to U$, and 
$\bar E \cap \bar f^{-1}(U)$ is its section. 
We may assume that $\bar S_1\subset \bar f^{-1}(U)$, so
$\bar f|_{\bar S_1}:\bar S_1\to \bar f(\bar S_1)$ is a (smooth) $\PP^1$-bundle.
Let $\ell$ be a fiber of the projection $\bar f|_{\bar S_1}:\bar S_1\to \bar f(\bar S_1)$
and let $\Sigma$ be its minimal section. Let $n:= -\Sigma^2$.
Restricting the first relation of \eqref{equation-q=3-K} to $\bar S_1$ we get
\[
-K_{\bar S_1}=(-K_{\bar X}-\bar S_1)|_{\bar S_1}=
2\bar S_1 |_{\bar S_1} +2\bar E|_{\bar S_1}=2\bar E|_{\bar S_1}+2\ell.
\]
Here $\bar E|_{\bar S_1}$ is an irreducible curve, a section of $\bar f|_{\bar S_1}:{\bar S_1}\to \bar f(\bar S_1)$.
On the other hand, 
\[
-n=-K_{\bar S_1}\cdot \Sigma-2= 2\bar E|_{\bar S_1}\cdot \Sigma+2\ell\cdot \Sigma-2=2\bar E|_{\bar S_1}\cdot \Sigma\ge 0.
\]
This gives us $n=0$, i.e. ${\bar S_1}\simeq \PP^1\times \PP^1$, and $\bar E|_{\bar S_1}\sim \Sigma$.
Further, the divisor $\bar S_2$ is $\bar f$-horizontal.
Indeed, otherwise $\bar S_2\sim 2\bar S_1$ and $\dim |\bar S_2|= \dim |\OOO_{\PP^2}(2)|=5$,
a contradiction. Since $a_2\ge 1$, the second relation of \eqref{equation-q=3-K}
can be written as $K_{\bar X}+\bar S_1+\bar S_2+ \bar E\sim 0$, so $\bar S_2\sim 2\bar S_1+\bar E$.
Clearly, $\bar E$ is a fixed component of the linear system $|\bar S_2-\bar S_1|= |\bar S_1+\bar E|$
(otherwise $\dim |\bar S_1+\bar E|>\dim |\bar S_1|=2$ and so $\dim |S_1|>2$).
Hence, $\dim |\bar S_2-\bar S_1|=\dim |\bar S_1|=2$.
From the exact sequence 
\[ 
0 \longrightarrow \OOO_{\bar X}(\bar S_2-\bar S_1)
 \longrightarrow \OOO_{\bar X}(\bar S_2)
 \longrightarrow \OOO_{\bar S_1}(\bar S_2)
 \longrightarrow 0
\]
we get 
\[
\dim H^0 ( \OOO_{\bar S_1}(\bar S_2))
\ge 
\dim H^0 (\OOO_{\bar X}(\bar S_2))- 
\dim H^0 ( \OOO_{\bar X}(\bar S_2-\bar S_1))
\ge 7.
\]
On the other hand, $\bar S_2|_{\bar S_1}\sim (2\bar S_1+\bar E)_{\bar S_1}\sim \Sigma+2\ell$ 
and so 
\[
\dim H^0 ( \OOO_{\PP^1\times \PP^1}(\Sigma+2\ell))=6, 
\]
a contradiction.
\end{scase}

\begin{scase}
Now let $\bar f$ be birational. Then
\[
 \hat q =3s_1+a_1e=2e+3(s_1+m_1e)\ge 5.
\]
On the other hand, $\g (\hat X)\ge 21$, so either $\hat q\le 7$ or
$\hat X\simeq \PP(1,2,3,5)$ by Theorem \ref{theorem-main-10-19}
and \ref{theorem-main-q=9}-\ref{theorem-main-q=6} of Theorem \ref{theorem-main}.
By Lemma \ref{lemma-torsion} the group $\Cl(\hat X)$ is torsion free.
\end{scase}

\begin{scase} \textbf{Subcase $\hat q\le 7$.}
Then $s_1=1$ and $e\le 2$. In particular, $\hat q\neq 6$.
By Proposition \ref{proposition-dim-A} we have $\hat X\simeq \PP(1^3,2)$.
Then $e=1$ and $a_1=2$.
Since $\dim |\Theta|=\dim |A|$, the linear system $|\Theta|$ is the 
proper transform of $|A|$.

As in the case \ref{case-q=3-933} using
\eqref{equation-discrepancies-second} we get 
\[
b=2+5\gamma_1=5\delta-3\ge 2. 
\]
Hence $Q:=\bar f(\bar F)$ is a smooth point.
Since $\Bs|\Theta|$ coincides with the point of index $2$,
we have 
$\bar f(\bar F)\not \subset \Bs|\Theta|$ and $\gamma_1=0$.
Then $b=2$ and $\delta=1$.
Let $\hat {\SSS}_1$ be the subsystem of $|\Theta|=|\hat S_1|$
consisting of divisors passing through $Q$.
Write
\[
\bar f^* \hat \SSS_1=\bar \SSS_1 +\gamma'_1\bar F,
\]
where $\bar \SSS_1\subset \bar X$ is the proper transform
and $\gamma_1'>0$.
Then 
\[
0\sim f^*(K_{\hat X}+3\hat \SSS_1 +2\hat E)\sim 
K_{\bar X}+3\bar \SSS_1 +2\bar E +3\gamma'_1\bar F.
\]
We get a contradiction with $q=3$.
\end{scase}

\begin{scase}
\textbf{Subcase $\hat X\simeq \PP(1,2,3,5)$.} 
Then $s_1\ge 3$ and $s_2\ge 7$.
From the second relation of \eqref{equation-q=3-K} we get 
$11=\hat q=s_1+s_2+e(1+m_1+m_2)$.
On the other hand, $9=\dim |7\Theta|\ge \dim |\bar S_2|$.
Thus $\dim |2A|=9$. This is possible only in the case \ref{no-q=3-B57}.
Moreover, $s_1=3$, $s_2= 7$, $e=1$, $a_1=2$, and $a_2=1$.

Further, $\dim |3\Theta|=\dim |A|$ and $|3\Theta|$ is the 
proper transform of $|A|$.
Similarly, $|7\Theta|$ is the 
proper transform of $|2A|$.
As in the case \ref{case-q=3-933} using
\eqref{equation-discrepancies-second} we get 
\[
b=\frac23 +\frac {11}3\gamma_1=\frac17 +\frac{11}7 \gamma_2\ge \frac 23. 
\]
Hence $Q:=\bar f(\bar F)$ is either a smooth point or a curve \cite{Kawamata-1996}
(otherwise $\bar f(\bar F)$ is a point of index $r=2$, $3$, or $5$, 
$\gamma_1\in \frac 1r\ZZ$ and $b>1/r$).
In this case $b$ is an integer, so $\gamma_1,\, \gamma_2>0$.
This means that $\bar f(\bar F)\subset \Bs|3\Theta|$ and $\bar f(\bar F)\subset \Bs|7\Theta|$.
On the other hand, $\Bs|7\Theta|\cap \Bs|3\Theta|$ is given by $\{x_1=x_3=x_2x_5=0\}\subset \PP(1,2,3,5)$.
Thus $\bar f(\bar F)$ is a point of index $2$ or $5$, a contradiction.
\end{scase}
\end{case}
Finally, as in previous sections, we need the following.

\begin{mtheorem}
\label{q=3-general-case}
{\bf Lemma.}
Let $X$ be a $\QQ$-Fano threefold with $\qQ(X)=3$.
If $\qW(X)\neq \qQ(X)$, then $\g(X)\le 16$. 
\end{mtheorem}

\begin{proof}
By \cite[Lemma 3.2 (iv)]{Prokhorov2008a} there is a 
$3$-torsion element $\Xi\in \Cl(X)$.
Let $P_1,\dots,P_k\in X$ be all the points where 
the divisor $\Xi$ is not Cartier and let $r_1,\dots r_k$ be 
indices of these points. According to \cite[Prop. 2.9]{Prokhorov2008a}
we have 
$\sum r_i=18$,
where each $r_i$ is divisible by $3$. Therefore, 
\[
-K_X\cdot c_2\le 24 -\sum (r_i-1/r_i) =6+\sum 1/r_i \le 8.
\]
By \cite[Proposition 2.2]{Suzuki-2004} 
\[
(4\qQ(X)-3) (-K_X)^3\le 4\qQ(X)^2 (-K_X\cdot c_2).
\]
Hence, $-K_X^3\le -4K_X\cdot c_2\le 32$.
By \cite[(2.5.3)]{Prokhorov-2007-Qe} we get
\[
\g(X)< \frac12 (-K_X)^3+2=17. 
\]
\end{proof}

\def\cprime{$'$} \def\mathbb#1{\mathbf#1} \def\bblapr{April}


\begin{thebibliography}{BKR10}

\bibitem[Ale94]{Alexeev-1994ge}
V. Alexeev.
\newblock General elephants of {${\bf Q}$}-{F}ano 3-folds.
\newblock {\em Compositio Math.}, 91(1):91--116, 1994.

\bibitem[BKR10]{Brown2010}
G.~{Brown}, M.~{Kerber}, and M.~{Reid}.
\newblock {Fano 3-folds in codimension 4, Tom and Jerry, Part I}.
\newblock {\em ArXiv e-prints}, September 2010.

\bibitem[GRD]{GRD}
G. Brown and others.
\newblock Graded ring database.
\newblock \url{http://grdb.lboro.ac.uk}.

\bibitem[BS07]{Brown-Suzuki-2007j}
G. Brown and K. Suzuki.
\newblock Computing certain {F}ano 3-folds.
\newblock {\em Japan J. Indust. Appl. Math.}, 24(3):241--250, 2007.

\bibitem[Fu86]{Furushima86}
 M. Furushima.
\newblock Singular del Pezzo surfaces and analytic compactifications of 
$3$-dimensional complex affine space $\mathbf{C}^3$.
\newblock {\em Nagoya Math. J.}, 104: 1--28, 1986.

\bibitem[IP99]{Iskovskikh-Prokhorov-1999}
V.~A. Iskovskikh and Yu. Prokhorov.
\newblock {\em Fano varieties. {A}lgebraic geometry. {V}.}, {\em
 Encyclopaedia Math. Sci.} \textbf{47},
\newblock Springer, Berlin, 1999.

\bibitem[Kaw92]{Kawamata-1992bF}
Y. Kawamata.
\newblock Boundedness of {$\mathbf{Q}$}-{F}ano threefolds.
\newblock In {\em Proceedings of the International Conference on Algebra, Part
 3 (Novosibirsk, 1989)}, {\em Contemp. Math.} \textbf{131}, 439--445,
 Providence, RI, 1992. Amer. Math. Soc.

\bibitem[Kaw96]{Kawamata-1996}
Y. Kawamata.
\newblock Divisorial contractions to {$3$}-dimensional terminal quotient
 singularities.
\newblock In {\em Higher-dimensional complex varieties (Trento, 1994)}, 
 241--246. de Gruyter, Berlin, 1996.
 
\bibitem[Kaw01]{Kawakita-2001}
M. Kawakita.
\newblock Divisorial contractions in dimension three which contract divisors to
 smooth points.
\newblock {\em Invent. Math.}, 145(1):105--119, 2001.

\bibitem[Kaw05]{Kawakita2005}
M. Kawakita.
\newblock Three-fold divisorial contractions to singularities of higher
 indices.
\newblock {\em Duke Math. J.}, 130(1):57--126, 2005.

\bibitem[Mor75]{Mori-1975}
S. Mori.
\newblock On a generalization of complete intersections.
\newblock {\em J. Math. Kyoto Univ.}, 15(3):619--646, 1975.

\bibitem[MP08]{Mori-Prokhorov-2008}
S. Mori and Y. Prokhorov.
\newblock On {$\mathbf Q$}-conic bundles.
\newblock {\em Publ. Res. Inst. Math. Sci.}, 44(2):315--369, 2008.

\bibitem[Pro]{Prokhorov-2010-qfano}
Y. Prokhorov.
\newblock Pari code for searching {${\bf Q}$}-{F}ano 3-folds.
\newblock \url{http://mech.math.msu.su/algebra/wiki/doku.php/staff:prokhorov:q-fano}.

\bibitem[Pro07]{Prokhorov-2007-Qe}
Y. Prokhorov.
\newblock The degree of {$\mathbf Q$}-{F}ano threefolds.
\newblock {\em Russian Acad. Sci. Sb. Math.}, 198(11):1683--1702, 2007.

\bibitem[Pro10]{Prokhorov2008a}
Y. Prokhorov.
\newblock {$\mathbf Q$}-{F}ano threefolds of large {F}ano index{,} {I}.
\newblock Doc. Math., J. DMV, 15, 843--872, 2010.


\bibitem[Rei]{Reid-graded-rings}
M.~Reid.
\newblock Graded rings and varieties in weighted projective space, preprint.

\bibitem[Rei87]{Reid-YPG1987}
M. Reid.
\newblock Young person's guide to canonical singularities.
\newblock In {\em Algebraic geometry, Bowdoin, 1985 (Brunswick, Maine, 1985)},
{\em Proc. Sympos. Pure Math.} \textbf{46}, 345--414. Amer. Math.
 Soc., Providence, RI, 1987.

\bibitem[RS03]{Reid2003}
M. Reid and K. Suzuki.
\newblock Cascades of projections from log del {P}ezzo surfaces.
\newblock In {\em Number theory and algebraic geometry}, {\em
 London Math. Soc. Lecture Note Ser.} \textbf{303}, 227--249. Cambridge Univ. Press,
 Cambridge, 2003.

\bibitem[San96]{Sano-1996}
T. Sano.
\newblock Classification of non-{G}orenstein {${\bf Q}$}-{F}ano {$d$}-folds of
 {F}ano index greater than {$d-2$}.
\newblock {\em Nagoya Math. J.}, 142:133--143, 1996.

\bibitem[Sho85]{Shokurov-1985}
V.~V. Shokurov.
\newblock A nonvanishing theorem.
\newblock {\em Izv. Akad. Nauk SSSR Ser. Mat.}, 49(3):635--651, 1985.

\bibitem[Suz04]{Suzuki-2004}
K. Suzuki.
\newblock On {F}ano indices of {$\mathbf{Q}$}-{F}ano 3-folds.
\newblock {\em Manuscripta Math.}, 114(2):229--246, 2004.

\end{thebibliography}
\end{document}